%% file: arxiv.tex
\documentclass[review,onefignum,onetabnum]{siamonline190516}

\usepackage{bm,bbm,stmaryrd}
\providecommand{\mathbold}[1]{\bm{#1}}
\newcommand{\mathbb}[1]{\mathbbm{#1}}
\newcommand{\mtx}[1]{\mathbold{#1}}
\newcommand{\Id}{\mathbf{I}}

\newcommand{\RR}{\mathbbm{R}}

\usepackage[notcite,notref,final]{showkeys}

\DeclareMathOperator{\rank}{rank}

\usepackage{graphicx}
\usepackage{epstopdf}
\usepackage{diagbox}
\usepackage{pict2e}

\newcommand{\tsp}{\top}
\newcommand{\econst}{\mathrm{e}}
\newcommand{\iunit}{\mathrm{i}}
\newcommand{\diff}{\mathrm{d}}

\usepackage{multirow}

\newcommand*{\norm}[1]{\left\|#1\right\|}

\newcommand*{\lowrank}[1]{\llbracket #1 \rrbracket}
\nolinenumbers

\input{ex_shared}

\usepackage{hyperref}

\ifpdf
\hypersetup{
  pdftitle={Streaming Kernel Analog Forecasting},
  pdfauthor={D. Giannakis, A. Henriksen, J. A. Tropp, and R. Ward}
  pdffunding={DG acknowledges support from NSF DMS 1854383 and ONR MURI N00014-19-1-242. AH was funded by AFOSR MURI FA9550-19-1-0005, NSF DMS 1952735.  JAT was supported by ONR N00014-18-1-2363 and NSF DMS 1952777.
  RW acknowledges support from AFOSR MURI FA9550-19-1-0005, NSF DMS 1952735, and NSF IFML 2019844.}
}
\fi

\begin{document}

\maketitle

\begin{abstract}
Kernel analog forecasting (KAF) is a powerful methodology for data-driven, non-parametric forecasting of dynamically generated time series data.  This approach has a rigorous foundation in Koopman operator theory and it produces good forecasts in practice, but it suffers from the heavy computational costs common to kernel methods.  This paper proposes a \textit{streaming} algorithm for KAF that only requires a single pass over the training data.  This algorithm dramatically reduces the costs of training and prediction without sacrificing forecasting skill.  Computational experiments demonstrate that the streaming KAF method can successfully forecast several classes of dynamical systems (periodic, quasi-periodic, and chaotic) in both data-scarce and data-rich regimes.  The overall methodology may have wider interest as a new template for streaming kernel regression.
\end{abstract}



\section{Introduction}

Forecasting problems are ubiquitous in physical science and engineering applications, including climate prediction \cite{PalmerHagedorn06}, navigation \cite{qu2009cooperative}, and medicine \cite{jackson2015applications}.  In these settings, we do not possess complete information about the state of the system, and we may not have full knowledge of the equations of motion.  Owing to our lack of omniscience, it is not possible to make predictions by integrating the current state forward in time.  Instead, we may acquire training data by observing some aspect of the system's evolution.  The goal is to build a compact model of the dynamics of this observable.  Given a new observation, the model should allow us to forecast the future trajectory from the initial condition.

Kernel analog forecasting (KAF)~\cite{2020alexander} offers a promising approach to this problem.  KAF is a data-driven, non-parametric forecasting technique that is best understood as a type of regularized kernel regression (\cref{sec:kaf}).  KAF emerged from recent efforts~\cite{BerryEtAl20} to translate Koopman operator theory into effective computational methodologies for forecasting (\cref{sec:koop}).  The approach belongs to a rapidly expanding literature~\cite{Mezic05,FroylandEtAl14,BruntonEtAl17,LiEtAl17,KlusEtAl18} on operator-theoretic techniques for low-order modeling of dynamical systems, including methods~\cite{BerryEtAl15,WilliamsEtAl15,Kawahara16,HamziOwhadi21} based on kernels.

KAF is mathematically rigorous, and it provides good-quality predictions for benchmark examples~\cite{2020alexander}.  Nevertheless, the straightforward implementation (``na\"ive KAF'') has several weaknesses.  First, na\"ive KAF requires multiple views of the training data, so it cannot operate in the ``streaming'' setting where we only see the training data once (\cref{sec:streaming}).
Second, the process of constructing the model is computationally expensive: to form the kernel matrix, the costs of arithmetic and storage are both \textit{quadratic} in the length of the training data.  
Third, the basic method must store all of the training data to make predictions, so the forecasting model is quite large.  Fourth, the arithmetic cost of a single forecast is \textit{linear} in the amount of training data.  These issues have limited the applicability of the KAF methodology.

In response to this challenge, we propose a novel \textit{streaming KAF} algorithm (\cref{sec: scalable kaf}).  Our approach depends on two prominent techniques from the field of randomized matrix computation~\cite{martinsson2020randomized}:
random Fourier features~\cite{2008Rahimi} for kernel approximation and the randomized Nystr{\"o}m method~\cite{2011Halko,Git13:Topics-Randomized,LLS+17:Algorithm-971,2017Tropp,martinsson2020randomized} for streaming PCA.  Overall, the streaming KAF method builds a model using time and storage \textit{linear} in the amount of training data, and it can make forecasts with time and storage that are \textit{independent} of the amount of training data.

Computational experiments (\cref{sec: experiments}) demonstrate that streaming KAF is a practical method for making predictions of two benchmark dynamical systems: Lorenz '63 (L63) \cite{Lorenz63} and two-level Lorenz '96 (L96) \cite{fatkullin2004computational}.  In particular, streaming KAF exhibits forecasting skill similar to na{\"i}ve KAF in a range of situations, including systems that are periodic, quasi-periodic, and chaotic.  At the same time, streaming KAF can operate in settings where na{\"i}ve KAF is prohibitively expensive, including cases where the observables are high-dimensional or the amount of training data is enormous.  In the data-rich setting, after just a few minutes of training time, streaming KAF can drive the forecasting error toward zero.  As a consequence, we believe that the streaming KAF algorithm has the potential to unlock the full potential of KAF as a forecasting methodology.

\begin{remark}[Prior work]
Although developed independently, our methodology is related to recent papers that apply random features to perform
streaming kernel principal component analysis~\cite{2016Ghashami,2018Ullah} and kernel ridge regression~\cite{AKM+17:Random-Fourier,2017Rudi}.
The details of our algorithm are somewhat different from these works, and we believe that our work yields a novel approach for \textit{streaming} kernel regression.  We have also studied a streaming KAF algorithm based on AdaOja~\cite{2019Henriksen},
an adaptive variant of Oja's algorithm, which is a competitive alternative to the Nystr{\"o}m method~\cite{2021HenriksenThesis}.
See~\cref{sec:related} for more discussion of related work.
\end{remark}

\subsection{Outline}

\Cref{sec:kaf} motivates the existing KAF procedure as a form of regularized kernel regression that is specifically designed for dynamical systems.
\Cref{sec: scalable kaf} describes how to develop a streaming implementation of KAF.  In particular, we discuss kernel approximation
via random Fourier features and the randomized Nystr{\"o}m method. 
Finally, \cref{sec: experiments} presents computational experiments which demonstrate that our methodology is
effective for two classical dynamical systems. 

\subsection{Notation} Throughout, we work in a real Euclidean space $\RR^d$ equipped with the $\ell_2$ norm $\norm{\cdot}$ and inner product $\langle \cdot, \cdot\rangle$.  The methodology and results should extend to the complex field $\mathbb{C}$. 
Matrices (such as $\mtx M \in \RR^{d \times n}$) are written as bold capitals,  vectors ($\mtx v \in \RR^d$) are written as bold lowercase, and scalars ($x \in \RR$) are written in plain lowercase.

Every matrix $\mtx M \in \RR^{d \times n}$ admits a compact singular value decomposition (SVD), a matrix factorization $\mtx M = \mtx U \mtx \Sigma \mtx V^\tsp $ with the following properties.  For $r \coloneqq \rank(\mtx M) \leq \min\{ d, n\}$, the left and right \textit{singular vector} matrices $\mtx U \in \RR^{d \times r}$ and $\mtx V \in \RR^{r \times n}$ have orthonormal columns.  The matrix $\mtx \Sigma = \diag(\sigma_1, \sigma_2, \dots, \sigma_r) \in \RR^{r \times r}$ is positive and diagonal, with its diagonal elements (the \textit{singular values}) arranged in decreasing order: $\sigma_1 \geq \sigma_2 \geq \dots \geq \sigma_r > \sigma_{r+1} \coloneqq 0$.  Singular values are uniquely determined, but singular vectors are not.

Given an SVD of the rank-$r$ matrix $\mtx{M}$, we can define the Moore--Penrose {pseudoinverse} $\mtx M^{\dagger} \coloneqq \mtx V \mtx \Sigma^{-1} \mtx U^\tsp $ where $\mtx{\Sigma}^{-1} \coloneqq \diag(\sigma_1^{-1}, \dots, \sigma_r^{-1}) \in \RR^{r \times r}$.  The pseudoinverse coincides with the matrix inverse for a full-rank, square matrix.

The operator norm $\norm{\mtx{M}} \coloneqq \sigma_1$ equals the largest singular value $\sigma_1$.  The Frobenius norm $\norm{ \mtx{M} }_{\mathrm{F}} \coloneqq ( \sum_{i=1}^r \sigma_i^2 )^{1/2}$ is the $\ell_2$ norm of the singular values.

For any rank parameter $\ell \leq r$, we can construct an $\ell$-truncated SVD $\lowrank{\mtx{M}}_\ell \coloneqq \mtx U \mtx \Sigma_\ell \mtx V^\tsp$ where $\mtx{\Sigma}_\ell \coloneqq  \diag(\sigma_1, \sigma_2, \dots, \sigma_\ell, 0, \dots, 0) \in \RR^{r \times r}$ retains only the leading $\ell$ singular values.  The matrix $\lowrank{\mtx{M}}_\ell$ is a best rank-$r$ approximation of $\mtx{M}$ with respect to both the operator norm and the Frobenius norm.  In the case $\sigma_{\ell} = \sigma_{\ell + 1}$, the truncated SVD $\lowrank{\mtx{M}}_{\ell}$ depends on the underlying choice of SVD, so this notation should be interpreted with care.

\section{Introduction to KAF} \label{sec:kaf}

Suppose we have access to snapshots of a discrete dynamical system as it evolves in time, and we would like to forecast its future values.
Let us begin with the most basic setting; we will discuss more general observation models in \cref{sec:other-observables}.

To formalize the problem, let ${\cal M} \subseteq \RR^d$ be a closed subset of a Euclidean space.  We call ${\cal M}$ the \textit{state space}. Let $F: {\cal M} \rightarrow {\cal M}$ be a mapping, called the \textit{flow map}.  Suppose that we observe an initial condition $\mtx x_0 \in {\cal M}$ as well as the (partial) trajectory $\mtx x_1, \mtx x_2, \dots, \mtx x_{n-1} \in {\cal M}$ obtained by iterating the flow map:
\begin{equation} \label{eqn:state-trajectory}
\mtx x_j = F(\mtx x_{j-1}) = F^j(\mtx{x}_0) \quad\text{for $j = 1, 2, \dots, n-1$.}
\end{equation}
In most settings, we do not actually know the flow map $F$.  Rather, the goal is to use information latent in the measured trajectory $(\mtx{x}_0, \dots, \mtx{x}_{n-1})$ to infer the dynamics.  Afterward, we are given a new initial condition $\mtx{y} \in \mathcal{M}$, and we are asked to forecast the future state $F^q(\mtx{y})$ of the system after $q$ time steps.

\subsection{Linear forecasting}

To motivate the KAF method, we first describe an earlier approach to the forecasting problem, based on \emph{linear inverse models} (LIMs) \cite{Penland89} and the closely related \textit{dynamic mode decomposition} (DMD) \cite{RowleyEtAl09,Schmid10,TuEtAl14,KlusEtAl18}.
Fix a forecasting horizon 
$q \in \mathbb{N}$.  We can arrange the observed trajectory $\mtx x_0, \mtx x_1, \dots, \mtx x_{n  -1}, \mtx x_n, \dots, \mtx x_{n + q - 1} \in \RR^{d}$ into a training data set that consists of input--response pairs:  $\{ (\mtx x_{j}, \mtx x_{j+q}) \}_{j=0}^{n-1}$.  Equivalently, consider the pair of matrices
\begin{align}
\label{input-output}
\mtx X &\coloneqq \begin{bmatrix} \mtx x_0 & \mtx x_1 & \dots & \mtx x_{n-1} \end{bmatrix} \in \RR^{d \times n}; \\
\quad \mtx X_{[+q]} &\coloneqq \begin{bmatrix} \mtx x_q & \mtx x_{q+1} & \dots & \mtx x_{n+q-1} \end{bmatrix} \in \RR^{d \times n}.
\end{align}
We can attempt to find the best linear model $\mtx{A} : \RR^d \to \RR^d$
for the dynamics by means of a least-squares fit: 
\begin{equation}
\label{eq:interpolate}
\mtx A  \in \underset{\mtx M \in \RR^{d \times d}}{\arg \min}\ \sum_{j=0}^{n-1} \norm{ \mtx M \mtx x_j - \mtx x_{j+q} }^2
	=  \underset{\mtx M \in \RR^{d \times d}}{\arg \min}\ \norm{ \mtx M \mtx X - \mtx X_{[+q]} }_{\mathrm{F}}^2.
\end{equation}
An optimal solution to this problem is the matrix
\begin{equation}
    \label{eqDMDA}
    \mtx A = \mtx X_{[+q]} \mtx X^{\dagger} \in \RR^{d \times d}.
\end{equation}
Suppose we are given a state $\mtx y \in \RR^d$ that serves as a new initial condition.
We can forecast the state $\mtx{y}_{[+q]} \coloneqq F^q(\mtx{y})$ after $q$ time steps
via the estimate $\mtx{y}_{[+q]} \approx \mtx{A} \mtx{y}$.  In other words, $\mtx{A}$ serves as a linear approximation to the iterated flow map $F^q$.  

Let us manipulate the linear model for the dynamics so that it takes a more suggestive form.
Recall that the pseudoinverse satisfies
$\mtx{X}^\dagger = (\mtx{X}^\tsp \mtx{X})^{\dagger} \mtx{X}^\tsp$.  Therefore, 
$$
\mtx{A} = \mtx{X}_{[+q]} (\mtx{X}^\tsp \mtx{X})^{\dagger} \mtx{X}^\tsp.
$$
Given a new initial condition $\mtx{y} \in \RR^d$, we obtain the linear forecast
\begin{equation} \label{eqn:linear-forecast}
\tilde{\mtx{y}}_{[+q]} \coloneqq \mtx{A} \mtx{y} = \mtx{X}_{[+q]} (\mtx{X}^\tsp \mtx{X})^{\dagger} (\mtx{X}^\tsp \mtx{y}) \in \RR^d.
\end{equation}
Observe that this computation can be formulated in terms of inner products between states.

\subsection{The kernel trick} 

Of course, dynamical systems of practical interest are highly nonlinear, so linear approximations are only valid over a short time horizon.  When one needs to process data with nonlinear structure, a general principle is to ``lift and linearize''.  That is, we apply a nonlinear map to transport the data to a high-dimensional space where it may have linear structure; we implement a linear fitting algorithm on the high-dimensional space; and then we project back down to the original domain to obtain a (nonlinear) low-dimensional model for the data. This approach gives rise to KAF, discussed below, as well as other data-driven analysis and forecasting techniques  \cite{WilliamsEtAl15,Kawahara16,KlusEtAl19,HamziOwhadi21}. 

Remarkably, this lifting technique can often be implemented without applying the nonlinear map explicitly.  Consider a method, such as \cref{eqn:linear-forecast}, that processes Euclidean data using the inner product as a measure of the similarity between data points.  The \textit{kernel trick} allows us to develop a nonlinear extension simply by replacing each inner product $\mtx{x}^\tsp \mtx{y}$ in the data space with a more general function $\kappa(\mtx{x}, \mtx{y})$, called a \textit{kernel}.

The kernel trick is justified by the Moore--Aronszajn theorem~\cite[section~2(4)]{Aro50:Theory-Reproducing}.
Let $\kappa : \RR^d \times \RR^d \to \RR$ be a symmetric, positive-definite function.  That is,
\begin{equation*} \label{eqn:psd-kernel}
\begin{bmatrix} \kappa(\mtx{v}_i, \mtx{v}_j) \end{bmatrix}_{i,j = 1}^{n}
\quad\text{is positive definite for every $n \in \mathbb N$ and $\mtx{v}_1, \dots, \mtx{v}_{n} \in \RR^d$.}
\end{equation*}
Then, the kernel function $\kappa$ coincides with the inner product on a Hilbert space $\cal{H}$.  More precisely, there is a nonlinear \textit{feature map} $\varphi : \RR^d \to \cal{H}$ with the property that $\kappa(\mtx{x}, \mtx{y}) = \langle \varphi(\mtx{x}), \varphi(\mtx{y}) \rangle_{\mathcal{H}}$ for all $\mtx{x}, \mtx{y} \in \RR^d$.  Implicitly, the feature map summarizes each data point $\mtx{x}$ by a long list $\varphi(\mtx{x}) \in {\cal H}$ of features, and the kernel computes the inner product between the feature vectors.  

One of the most popular kernel functions is the Gaussian radial basis function (RBF) kernel.
For an inverse bandwidth parameter $\gamma > 0$,
this kernel takes the form
\begin{equation} \label{eqn:gauss-rbf}
\kappa(\mtx{x}, \mtx{y}) \coloneqq \econst^{- \gamma \norm{ \mtx{x} - \mtx{y} }^2}
\quad\text{for $\mtx{x}, \mtx{y} \in \RR^d$.}
\end{equation}
Under this kernel, two points are ``similar'' precisely when they are close enough together in Euclidean distance,
where the scale depends on the choice of $\gamma$.
For clarity of presentation, we will work exclusively with the Gaussian RBF kernel in this paper.

\subsection{Nonlinear kernel forecasting}

We can apply the kernel trick to the linear forecasting model~\cref{eqn:linear-forecast}.
Indeed, we may replace the inner-products in the forms $\mtx{X}^\tsp \mtx{X}$ and $\mtx{X}^\tsp \mtx{y}$
by their kernel equivalents:
$$
\mtx{K}_{x,x} \coloneqq \big[ \kappa(\mtx{x}_i, \mtx{x}_j) \big]_{i,j} \in \RR^{n \times n}
\quad\text{and}\quad
\mtx{K}_{x,y} \coloneqq \big[ \kappa(\mtx{x}_i, \mtx{y}) \big]_{i} \in \RR^{n}.
$$
This step leads to the kernel analog forecast
\begin{equation} \label{kernel_full}
f_q(\mtx{y}) \coloneqq \mtx{X}_{[+q]} (\mtx{K}_{x,x})^\dagger \mtx{K}_{x,y} \in \RR^d.
\end{equation}
The forecast~\cref{kernel_full} provides a natural nonlinear
generalization of the linear forecast \cref{eqn:linear-forecast}.

\subsection{Regularization}

It is dangerous to implement the formula~\cref{kernel_full}
as written because kernel matrices, such as $\mtx{K}_{x,x}$,
are notoriously ill-conditioned; for example, see~\cite{belkin2018approximation}.
As a consequence, the method~\cref{kernel_full} can be sensitive to small changes
in the observed data.

The paper~\cite{2020alexander} proposes a mechanism for stabilizing
the nonlinear forecast~\cref{kernel_full} by replacing the kernel matrix $\mtx{K}_{x,x}$ with
its best rank-$\ell$ approximation $\lowrank{\mtx{K}_{x,x}}_{\ell}$,
where $\ell \in \mathbb{N}$ is a parameter. 
In practice, we must also shift the kernel matrix by $\mu \Id$ by
a small parameter $\mu$ to avoid numerical problems.
These modifications leads to the stabilized kernel analog forecast
\begin{equation}
    \label{kernel_full_2}
f_{q,\ell}(\mtx Y) \coloneqq  \mtx X_{[+q]}  (\lowrank{\mtx K_{x,x} + \mu \Id}_{\ell})^{\dagger} \mtx K_{x,y}.
\end{equation}
The dimension $\ell$ of the regression model is usually modest (say, 100s or 1000s);
it increases slowly with the required accuracy of the forecasts.
The shift parameter $\mu$ is taken to be a small fixed value, such as $10^{-6} \norm{ \mtx{K}_{x,x} }$.

The forecasting method~\cref{kernel_full_2} is rigorously justified in~\cite{2020alexander}.
We can view the approach as a form of regularized least-squares~\cite{hanson1971numerical, varah1973numerical, bjorck1996numerical}
on the feature space induced by the kernel.  It is closely related to
kernel ridge regression~\cite{scholkopf2002learning}.

\subsection{Resource usage}
\label{sec:kaf-cost}

The KAF method~\cref{kernel_full_2} involves two phases.
In the training step, we use the trajectory data $\mtx{X}$ to compute
a matrix of prediction weights.  In the forecasting step,
we use the trajectory data and the test state $\mtx{y}$
to make the forecast.
Let us summarize the resource usage of an uninspired implementation
of the KAF procedure (``na{\"i}ve KAF'').  See \cref{tab:kaf} for a summary of this discussion.

In the training phase, we first construct the $n \times n$
kernel matrix $\mtx{K}_{x,x}$.  This step involves $O( dn^2 )$
arithmetic and $O(n^2)$ storage.  The quadratic dependency on the number $n$
of training samples is a severe bottleneck that prevents us from performing KAF at scale.

Next, we must compute the $\ell$-truncated eigenvalue decomposition of the kernel matrix $\mtx{K}_{x,x}$.
Classical algorithms can succeed with $O(\ell^2 n)$ arithmetic operations
and $O(\ell n)$ storage.  Nevertheless, dense methods require random access to
the kernel matrix, while Krylov methods require a long sequence of
matrix--vector multiplies with the kernel matrix~\cite{GvL13:Matrix-Computations}.
Moreover, these algorithms are not fully reliable~\cite{LLS+17:Algorithm-971}.

Third, we form the matrix $\mtx{W} \coloneqq \mtx{X}_{[+q]} (\lowrank{ \mtx{K}_{x,x} + \mu \Id }_{\ell})^\dagger \in \RR^{d \times n}$
of prediction weights.  Using the factorized form of the eigenvalue decomposition,
this product costs $O(d\ell n)$ operations.  The weight matrix
requires storage $O(d n)$, which is comparable to the cost of
storing the original trajectory data.

To make a forecast from a single initial condition $\mtx{y} \in \RR^d$,
we need to perform the kernel computation $\mtx{K}_{x,y} \in \RR^n$.
The cost is $O(dn)$ operations and $O(n)$ storage.
To complete the forecast, we form the matrix--vector product
$\mtx{W}\mtx{K}_{x,y}$,
at a cost of $O(dn)$ operations.
The linear dependency on the number $n$ of training
points means that forecasting is very expensive.

\subsection{Other observables}
\label{sec:other-observables}

The KAF methodology extends to a wider setting.
\Cref{sec: scalable kaf} provides full details for a streaming KAF
algorithm at this level of generality. For now, we just sketch the idea.

Suppose that we observe the value of
a function $u : \mathcal{M} \to \mathcal{N}$ of the state, which is called a \textit{covariate}.
For simplicity, we will always take $\mathcal{N} = \RR^{d'}$.
Given an observed covariate $u(\mtx{x})$, we would like to predict a function $g : \mathcal{M} \to \RR^r$
of the state $\mtx{x}$, which is called a \textit{response variable}.  Functions of the state,
such as $g$ and $u$, are called \textit{observables}.%
\footnote{It is important the the response variable $g$ takes values in a linear space.
In principle, the covariates $u$ could take values in a nonlinear manifold $\mathcal{N}$, but
we will not consider this extension.}

We can build a kernel analog forecast for future values of the response
by introducing a kernel $\tilde{\kappa} : \RR^{d'} \times \RR^{d'}$ on the covariate space.
Roughly speaking, we replace the matrix $\mtx{X}$ of training state data by
observed covariate values $[ u(\mtx{x}_0), ..., u(\mtx{x}_{n-1}) ] \in \RR^{d' \times n}$.
Replace the matrix $\mtx{X}_{[q]}$ of lagged state data by
the lagged matrix $[ g(\mtx{x}_q), ..., g(\mtx{x}_{n+q-1}) ] \in \RR^{r\times n}$
of observed response variables.  Repeat the derivation above to obtain a KAF function $g_{\ell, q}$
for predicting the observable $g$ from the covariate $u$.

The computational costs are similar to the costs of the basic KAF method,
but the state dimension $d$ is replaced by either the covariate dimension $d'$
or the response variable dimension $r$, depending on the role of the state in the computation.
See Table~\ref{sec:kaf} for an accounting.

\subsection{Connection with Koopman operator theory}
\label{sec:koop}

The linear approach~\cref{eqn:linear-forecast} to forecasting was originally proposed
in the paper \cite{Penland89}, and
the nonlinear kernel forecast~\cref{kernel_full}
was presented in~\cite{williams2015kernel}.
The paper \cite{TuEtAl14} clarifies the connection between the nonlinear forecast
and Koopman operator theory \cite{EisnerEtAl15}. The paper \cite{2020alexander} shows that KAF approximates the expectation of the response variable under the action of the Koopman operator, conditioned on the covariate data observed at forecast initialization.
Here is an informal summary of these ideas.

In plain language, the classical work of Koopman and von Neumann \cite{koopman1931hamiltonian,KoopmanVonNeumann32} characterizes a dynamical system through its induced action on a \emph{linear} space of observables. As a basic example, a real-valued function
$g : {\cal M} \to \RR$ on the state space is an observable 
of the dynamical system.  The Koopman operator ${\cal K}$ is a linear operator
on the space of observables that acts by composition with the flow map of the dynamics: $( {\cal K} g )(\mtx{x}) \coloneqq (g \circ F)(\mtx{x}) = g(F(\mtx{x}))$.
Regardless of the complexity of the dynamical system, we can understand its behavior
by spectral analysis of the linear operator ${\cal K}$ on an appropriately chosen Banach space of observables \cite{Baladi00,EisnerEtAl15}. In particular, since our state space $\mathcal M $ is a subset of $\mathbb R^d$, we can represent every state $ \mtx x \in \mathcal M $ by the ``identity'' observable, $\iota : \mathcal M \to \mathbb R^d $ with $\iota(\mtx{x}) = \mtx{x}$. Thus, the dynamical system becomes linear when lifted to a sufficiently high-dimensional space of observables: $ F( \mtx x  ) = (\mathcal K \iota)(\mtx x)$.  Using similar ideas, we can also represent dynamical systems with infinite-dimensional state spaces by means of linear Koopman operators.  

Building on previous work \cite{zhao2016analog, alexander2017kernel, comeau2017data},
the recent paper \cite{2020alexander} established that the stabilized forecast~\cref{kernel_full_2}
is a rigorous approximation of the Koopman dynamics of observables in the limit of large data.
Consider a measure-preserving and ergodic dynamical system $F$,
and let $[\mtx{x}_0, \dots, \mtx{x}_{n-1}]$ be a state trajectory as in~\cref{eqn:state-trajectory}.
Suppose we acquire training data in the form of covariate--response pairs
$(\mtx{u}_0, \mtx{g}_q), \dots, (\mtx{u}_{n-1}, \mtx{g}_{q+n-1})$,
where $\mtx{u}_i = u(\mtx{x}_i) \in \mathcal{N}$ and $\mtx{g}_i = g(\mtx{x}_i) \in \RR^r$.
We may construct the KAF function $g_{\ell, q}$ as summarized in~\cref{sec:other-observables}.
Let $\mtx{y} \in \mathcal{M}$ be an initial condition with an observed covariate $u(\mtx{y})$.
Then the kernel analog forecast converges%
\footnote{Convergence takes place in the $L_2$ norm of the invariant measure in the iterated limit of $\ell \to \infty$ after $n \to \infty$, and almost surely with respect to the initial condition $\mtx x_0$ in the training data.}
to the conditional expectation of the response under the Koopman operator,
given the covariate data at forecast initialization:
$$
g_{q, \ell}(\mtx{v}) \to \mathbb{E}[ (\mathcal{K}^q g)(\mtx{y}) \, \vert \, u(\mtx{y}) = \mtx{v} ]
\quad\text{in $L_2$ as $\ell, n \to \infty$.}
$$
The conditional expectation is the optimal $L_2$ approximation
to the Koopman evolution $(\mathcal{K}^q g)(\mtx{y})$,
given only the measured covariate $\mtx{v}$.
In the specific case where the observables $g = u = \iota $ reproduce the full state vector,
we deduce that the forecast $f_{q,\ell}(\mtx y)$ presented in~\cref{kernel_full_2}
converges to the true $q$-step dynamical evolution, $F^q(\mtx y)$.

\section{Streaming KAF}
\label{sec: scalable kaf}

While KAF is rigorously justified in the limit of large data, it also becomes prohibitively
expensive to implement because of its storage and arithmetic costs (\cref{sec:kaf-cost}).  Indeed, the time
required to construct the kernel matrix is \textit{quadratic} in the length $n$ of training data.
The time required to make a single forecast is \textit{linear} in $n$.
Furthermore, we need multiple views of the training data to build the
model and another view to make a forecast,
so the algorithm cannot operate in the streaming setting.

In this section, we will develop a streaming KAF method that resolves each of these issues.
Our algorithm processes the trajectory data in a single pass.
It reduces the arithmetic cost of training to be \textit{linear} in the number $n$
of training points, and the cost of each forecast becomes \textit{independent} of the
amount of training data.  It also limits the storage needed for the computations
and for the forecasting model.  Experiments (\cref{sec: experiments}) show that the
streaming KAF method is competitive with the original KAF method in forecasting skill
on problem sizes where the original KAF method is tractable.  But streaming KAF can
achieve significantly \textit{better} forecasts than na{\"i}ve KAF because the streaming
method can ingest large amounts of training data and resolve the dynamics more accurately.

\subsection{Streaming data}
\label{sec:streaming}

Streaming data models have become popular for working with time series that have many elements,
especially in high dimensions or in cases where the data arrives
at high velocity~\cite{muthukrishnan2005data}. 
The key features of a streaming data model%
\footnote{More general streaming models describe a sequence of \textit{update operations}
to a data domain.}
are that
(1) the elements of the time series are presented in sequential order;
(2) we must process each datum at the time it arrives; and
(3) we do not have sufficient storage to maintain the entire time series.
The goal is to extract enough information to answer a particular set of
questions about the observed data.
These constraints necessitate algorithms that can handle each element
individually and that build a compact representation of the time series
to support subsequent queries.

Streaming models are well suited to dynamical systems data
that has an explicit temporal order.
It would be appealing to scan linearly
through the trajectory data $(\mtx{x}_0, \mtx{x}_1, \dots, \mtx{x}_{n-1})$
a single time, discarding each state after we have processed it.
Our aim is to build a forecasting model that can take a query
state and predict the subsequent trajectory of the system.
Ideally, the forecasting model should be much smaller than the original training data.
Yet the basic KAF method fails this desideratum.  We will show how to accomplish this task.

\subsection{Overview}

Our streaming KAF method is based on two techniques from the field
of randomized matrix computations~\cite{martinsson2020randomized}.
First, we use random Fourier features (RFF) to build a
structured approximation of the original kernel function.
This approximation allows us to rewrite the KAF target
function~\cref{kernel_full}, replacing the $n \times n$
kernel matrix $\mtx{K}_{x,x}$ by a much smaller
matrix that is easier to compute and captures the same information.
This reformulation also allows us to avoid the kernel computation $\mtx{K}_{x,y}$,
which couples the training and test data.
As a consequence, we can build a more compact forecasting model.

When we restructure the KAF target function, the low-rank
approximation of the kernel matrix converts into a low-rank
approximation of the covariance matrix of the features of the training data.
The latter approximation may be interpreted as a streaming PCA problem.
Here, we employ the randomized Nystr{\"o}m method devised by
Halko et al.~\cite{2011Halko,Git13:Topics-Randomized,LLS+17:Algorithm-971}
and extended to the streaming setting in~\cite{2017Tropp,martinsson2020randomized}.
This algorithm requires minimal storage and arithmetic,
and it reliably produces a more accurate solution than competing methods.

The rest of this section introduces the random features construction.
It shows how to integrate random features into KAF to obtain a
streaming algorithm, and it highlights the role of the Nystr{\"o}m method.
Last, we compare the resource usage of streaming KAF with the direct implementation of KAF.
See~\cref{sec:related} for related work.

\subsection{Kernel approximation by random features}
\label{sec:rff}

Random Fourier features (RFF)~\cite{2008Rahimi} offer a simple and effective way to approximate
certain types of kernels, including the Gaussian RBF kernel.  This section summarizes the RFF
construction, and the next section explains how we can use RFF to forecast a dynamical system.

Bochner's theorem~\cite{1932Bochner} provides the mathematical foundation for RFF.
Let us consider a bounded, continuous, positive-definite kernel $\kappa: \RR^d \times \RR^d \to \RR$ on a Euclidean space.
Assume that the kernel is also translation invariant: $\kappa(\mtx{x}, \mtx{y}) \coloneqq h( \mtx{x} - \mtx{y} )$.
The theorem asserts that the kernel is the Fourier transform of a bounded positive measure.  More precisely,
there exists a unique probability measure $\nu$ on $\RR^d$ and a positive constant $c \coloneqq h({\bf 0})$ for which
$$
\kappa(\mtx{x}, \mtx{y}) = \int_{\RR^d}  c \, \diff \nu(\mtx{z})\, \econst^{\iunit \, \mtx{z}^\tsp (\mtx{x} - \mtx{y})}
	= \int_{\RR^d} c \, \diff \nu(\mtx{z}) \, (\econst^{\iunit \, \mtx{z}^\tsp \mtx{x}})( \econst^{\iunit \, \mtx{z}^\tsp \mtx{y}} )^*,
$$
where ${}^*$ denotes the complex conjugate.  Since we are working in the real setting, we can rewrite the last expression to avoid
complex-valued functions:
$$
\kappa(\mtx{x}, \mtx{y}) = \int_{\RR^d} 2c\, \diff \nu(\mtx{z}) \int_0^{2\pi} \frac{\diff{\theta}}{2 \pi} \,
	\cos( \theta + \mtx{z}^\tsp \mtx{x} ) \, \cos( \theta + \mtx{z}^\tsp \mtx{y}).
$$
This statement follows by direct calculation using trigonometric identities.  The key property of these
formulas is that the integrand is a \textit{separable} function of the variables $\mtx{x}$ and $\mtx{y}$.

The simple idea behind RFF is to approximate the kernel using a Monte Carlo estimate of the integral.
Let the parameter $s \in \mathbb{N}$ designate the number of random features.
Once and for all, draw and fix
independent random vectors $\mtx{z}_1, \dots, \mtx{z}_s \in \RR^d$ that are
distributed according to the probability measure $\nu$.
Draw and fix independent random scalars $\theta_1, \dots, \theta_s \in \RR$
with the $\textsc{uniform}[0,2\pi)$ distribution.  Then we can construct
a separable, rank-$s$ approximation $\hat{\kappa} : \RR^d \times \RR^d \to \RR$
of the original kernel:
$$
\hat{\kappa}(\mtx{x}, \mtx{y}) \coloneqq \frac{2c}{s} \sum_{i=1}^s \cos(\theta_i + \mtx{z}_i^\tsp \mtx{x}) \, \cos(\theta_i + \mtx{z}_i^\tsp \mtx{y} ).
$$
It is not hard to see that $\hat{\kappa}(\mtx{x}, \mtx{y}) \approx \kappa(\mtx{x}, \mtx{y})$ with high probability
for a fixed pair $(\mtx{x}, \mtx{y})$.

Equivalently, we may define a feature map $\varphi : \RR^d \to \RR^s$ by the formula
$$
\varphi(\mtx{x}) \coloneqq \sqrt{\frac{2c}{s}} \cdot  \big[ \cos(\theta_i + \mtx{z}_i^\tsp \mtx{x}) \big]_{i=1}^s.
$$
Then we can compute the approximate kernel $\hat{\kappa}$ as the inner product between two feature vectors:
$$
\hat{\kappa}(\mtx{x}, \mtx{y}) = \varphi(\mtx{x})^\tsp \varphi(\mtx{y}).
$$
In other words, the approximate kernel is a bilinear function of nonlinear features.

In computational settings, we are usually interested in approximating the kernel matrix $\mtx{K}_{x,x}$
associated with a family $\{ \mtx{x}_0, \dots, \mtx{x}_{n-1} \} \subset \RR^d$ of data points.  That is,
$$
\mtx{K}_{x,x} \coloneqq \big[ \kappa( \mtx{x}_i, \mtx{x}_j ) \big]_{i,j}
	\quad\approx\quad \big[ \hat{\kappa}( \mtx{x}_i, \mtx{x}_j ) \big]_{i,j} \eqqcolon \hat{\mtx{K}}_{x,x}.
$$
To this end, we collect the data points as the columns of a matrix $\mtx{X} \in \RR^{d \times n}$.
Extend the feature map $\varphi$ to matrices by applying the vector feature map to each \textit{column}.
Thus, $\varphi : \RR^{d \times n} \to \RR^{s \times n}$.  With this notation, we find that
$$
\hat{\mtx{K}}_{x,x} = \varphi(\mtx{X})^\tsp \varphi(\mtx{X}).
$$
The kernel matrix approximation is the Gram matrix of the nonlinear features.

Finally, we must discuss the number $s$ of random features that we need to ensure that the
kernel matrix approximation $\hat{\mtx{K}}_{x,x}$ serves in place of the true kernel matrix $\mtx{K}_{x,x}$
for machine learning tasks.  When we have $n$ training points, it has been shown~\cite{sriperumbudur2015optimal, 2017Rudi, 2018Ullah, szabo2019kernel} 
that it suffices to use
\begin{equation}
    \label{eq:rfeatures}
\text{$s = O( \sqrt{n} \, \log(n))$ random features}
\end{equation}
for kernel principal component analysis (KPCA) or for kernel ridge regression (KRR).
The justification involves statistical assumptions on the training and test data. 
Our empirical study indicates that, in our application, we may extract even fewer
features without much loss in forecasting performance.

As a particular example of the RFF construction,
consider the Gaussian RBF kernel~\cref{eqn:gauss-rbf} on $\RR^d$ with inverse bandwidth $\gamma > 0$.
The normalization constant $c = 1$, and the associated spectral measure $\nu$ satisfies
$$
\diff \nu(\mtx{z}) = (4 \pi \gamma)^{-d/2} \econst^{- \norm{\mtx{z}}^2 / (4\gamma)} \, \diff{\mtx{z}}.
$$
That is, the random feature descriptor $\mtx{z}$ is a centered normal vector with
covariance $(2\gamma) \Id$.

Algorithm~\ref{alg:rff} contains basic pseudocode for implementing Gaussian RBF random features. 
In this version, the feature descriptors require $O(ds)$ storage, and it costs $O(ds)$ operations
to compute the features for a single input vector.
The pseudocode also includes several methods for streaming computation of matrix--matrix products
with featurized data $\varphi(\mtx{X})$.

\begin{remark}[More efficient Gaussian feature maps] \label{rem:fast-features}
We can accurately approximate the RFF map for the Gaussian RBF kernel
using randomized trigonometric transforms~\cite{le2013fastfood, CN2021}.
This construction reduces the storage cost for the random feature descriptors to $O(s)$,
and it costs $O(s \log d)$ operations to compute the features of a single input vector.
For high-dimensional state spaces (or covariates), we can
obtain significant gains, but the basic construction is
superior in low-dimensional settings.
\end{remark}

\begin{remark}[Kernels that admit random feature maps]
It is also possible to construct random features for other
kinds of kernel functions, including kernels that are not translation invariant.
See~\cite[Sec.~19]{martinsson2020randomized} for some discussion and references.
\end{remark}

\begin{algorithm}[t!] 
  \caption{\textsl{Random Fourier Features for Gaussian RBF Kernel.}
  See~\cref{sec:rff}.}
  \label{alg:rff}
  
  \vspace{0.5pc}

  The constructor (RFF) generates a random feature map $\varphi$
  for the Gaussian RBF kernel on $\RR^d$ with inverse bandwidth $\gamma > 0$ with $s$ random features.
  The \textsc{Featurize} method of $\varphi$ applies the random feature map to the columns
  of the input matrix $\mtx{X} \in \RR^{d \times B}$ to obtain
  $\varphi(\mtx{X}) \in \RR^{s \times B}$.
  The other methods featurize an input matrix $\mtx{X} \in \RR^{d \times B}$
  and compute various matrix products between $\varphi(\mtx{X})$ and
  another input $\mtx{M}$ by streaming columns of $\mtx{X}$.
  
  \begin{algorithmic}[1]

\vspace{1pc}

  	\State \textbf{local variables} $\gamma \in \RR_{++}$ and $d, s \in \mathbb{N}$
		\Comment	RFF parameters
	\State \textbf{local variables} $\mtx{z}_1, \dots, \mtx{z}_s \in \RR^d$ and $\theta_1, \dots, \theta_s \in \RR$
		\Comment	Feature descriptors

  \vspace{0.5pc}

	\Function{RFF}{$\gamma \in \RR_{++}, d \in \mathbb{N}; s \in \mathbb{N}$}
		\Comment	{Initialization}
    \State	Store RFF parameters $\gamma, d; s$
    \For {$i = 1, \dots, s$}
    \State		$\mtx{z}_i \gets \sqrt{2 \gamma} \cdot \texttt{randn}(d, 1)$
    	\Comment	{Draw Gaussian vector}
    \State		$\theta_i \gets 2\pi \cdot \texttt{rand}(1,1)$
    	\Comment	{Draw uniform scalar}
    \EndFor
    \State	\Return \texttt{self}
    	\Comment	{Return feature map}
	\EndFunction

  \vspace{0.5pc}
	
	\Function{Featurize}{$\mtx{X} \in \RR^{d \times B}$}
		\Comment	{Compute features of $\mtx{X}$}
		
		\For {$j = 1, \dots, B$}
		\For {$i = 1, \dots, s$}
		\State	$[\varphi(\mtx{X})]_{ij} \gets \sqrt{2 / s} \cdot \cos(\theta_i + \mtx{z}_i^\tsp \mtx{X}(:,j) )$
		\EndFor
		\EndFor
		
		\State \Return $\varphi(\mtx{X}) \in \RR^{s \times B}$
	\EndFunction
	
	  \vspace{0.5pc}
	
	\Function{MultCov}{$\mtx{X} \in \RR^{d\times B}, \mtx{M} \in \RR^{s \times \ell}$}
		\Comment	{Form product $\varphi(\mtx{X}) \varphi(\mtx{X})^\tsp \mtx{M}$}
		
		\State	$\mtx{T} \gets \texttt{zeros}( s, \ell )$
		\For {$j = 1, \dots, B$}
						\Comment	Block for efficiency
		\State	$\mtx{v} \gets \textsc{Featurize}(\mtx{X}(:, j))$
			\Comment Compute features
		\State	$\mtx{T} \gets \mtx{T} + \mtx{v} \, (\mtx{v}^\tsp \, \mtx{M})$
		\EndFor
		
		\State \Return $\mtx{T} \in \RR^{s \times \ell}$
	\EndFunction

	\vspace{0.5pc}
	
		\Function{RMultAdj}{$\mtx{X} \in \RR^{d \times B}, \mtx{M} \in \RR^{r \times B}$}
		\Comment	{Form product $\mtx{M} \, \varphi(\mtx{X})^\tsp$}
		
		\State	$\mtx{T} \gets \texttt{zeros}( r, s)$
		\For {$j = 1, \dots, B$}
						\Comment	Block for efficiency
		\State	$\mtx{T} \gets \mtx{T} + \mtx{M}(:, j) \, \textsc{Featurize}(\mtx{X}(:, j))^\tsp$
		\EndFor
		
		\State \Return $\mtx{T} \in \RR^{R \times s}$
	\EndFunction

	\vspace{0.5pc}
	
		\Function{RMult}{$\mtx{X} \in \RR^{d \times B}, \mtx{M} \in \RR^{r \times s}$}
		\Comment	{Form product $\mtx{M} \, \varphi(\mtx{X})$}
		
		\State	$\mtx{T} \gets \texttt{zeros}( r, B )$
		\For {$j = 1, \dots, B$}
						\Comment	Block for efficiency
		\State	$\mtx{T}(:, j) \gets \mtx{M} \cdot \textsc{Featurize}(\mtx{X}(:, j))$
		\EndFor
		
		\State \Return $\mtx{T} \in \RR^{r \times B}$
	\EndFunction

\vspace{0.5pc}
\end{algorithmic}

\end{algorithm}

\subsection{KAF with random features}

We can use RFF to approximate the kernel matrices that appear in the
regularized KAF target function~\cref{kernel_full_2}.  Recall that the matrix $\mtx{X} \in \RR^{d\times n}$
contains the training data, while $\mtx{y} \in \RR^{d}$ is a piece of test data.
Draw and fix a random feature map $\varphi : \RR^d \to \RR^s$ with $s$ random features.
Then we can approximate the KAF as
\begin{equation} \label{eqn:kaf-rff}
\begin{aligned}
f_{q,\ell}(\mtx{y}) \quad\approx\quad \hat{f}_{q, \ell}(\mtx{y})
&\coloneqq \mtx{X}_{[+q]} ( \lowrank{ \hat{\mtx{K}}_{x,x} + \mu \Id }_{\ell} )^{\dagger} \hat{\mtx{K}}_{x,y} \\
&= \mtx{X}_{[+q]} ( \lowrank{ \varphi(\mtx{X})^\tsp \varphi(\mtx{X})  + \mu \Id }_{\ell} )^{\dagger} \varphi(\mtx{X})^\tsp \varphi(\mtx{y}) \\
&\eqqcolon \mtx{W}_{q, \ell} \cdot \varphi(\mtx{y}).
\end{aligned}
\end{equation}
The forecasting model consists of the matrix $\mtx{W}_{q, \ell} \in \RR^{d\times s}$ of prediction weights, along with the
description of the feature map $\varphi : \RR^{d} \to \RR^s$. 
A key benefit of the reformulation~\cref{eqn:kaf-rff} is the
complete decoupling of the test data $\mtx{y}$ from the forecasting model.

Direct substitution of random features does not lead immediately to a streaming algorithm.
Indeed, the formula~\cref{eqn:kaf-rff} involves the rank truncation of the
$n \times n$ approximate kernel matrix $\varphi(\mtx{X})^\tsp \varphi(\mtx{X})$.
We cannot form this matrix without multiple views of the columns of $\mtx{X}$,
and the matrix imposes unacceptable storage and arithmetic costs.

\subsection{Streaming KAF}

To develop a streaming algorithm, we first recast the expression~\cref{eqn:kaf-rff}
in terms of a much smaller $s \times s$ matrix. Recall the linear-algebraic identity
$$
(\lowrank{\mtx{M}^\tsp \mtx{M} + \mu \Id}_\ell)^\dagger \mtx{M}^\tsp
	= \mtx{M}^\tsp (\lowrank{\mtx{M}\mtx{M}^\tsp + \mu \Id}_{\ell})^\dagger.
$$
Using this formula, we can write the prediction weights as
\begin{equation} \label{eqn:scalable-weights}
\mtx{W}_{q,\ell} = (\mtx{X}_{[+q]} \, \varphi(\mtx{X})^\tsp) ( \lowrank{\varphi(\mtx{X}) \varphi(\mtx{X})^\tsp + \mu \Id}_{\ell} )^\dagger.
\end{equation}
The matrices in parentheses have the dimensions $d \times s$ and $s \times s$, respectively.
Moreover, this representation now supports a streaming algorithm.

In sequence, we pass over the columns $\mtx{x}_i$ of the training states,
generating random features $\varphi(\mtx{x}_i)$ on the fly.
Simultaneously, we update the covariance of the features and the covariance between the
features and the lagged data.  Beginning with $\mtx{C}_{xx} = \mtx{0}_{s \times s}$
and $\mtx{C}_{gx} = \mtx{0}_{d\times s}$, iterate
\begin{equation} \label{eqn:covariance-stream}
\mtx{C}_{xx} \gets \mtx{C}_{xx} + \varphi(\mtx{x}_i) \varphi(\mtx{x}_i)^\tsp
\quad\text{and}\quad
\mtx{C}_{gx} \gets \mtx{C}_{gx} + \mtx{x}_{i+q} \varphi(\mtx{x}_i)^\tsp.
\end{equation}
[Because of the lag, to form the matrix $\mtx{C}_{gx}$,
the algorithm must buffer the input states at a cost of $O(qd)$.]
Once we have streamed all of the training data, we may construct the matrix
of prediction weights as
\begin{equation} \label{eqn:streaming-weights}
\mtx{W}_{q, \ell} = \mtx{C}_{gx} \cdot (\lowrank{ \mtx{C}_{xx} + \mu \Id }_{\ell})^\dagger. 
\end{equation}
Since the expressions for the weights
in~\cref{eqn:kaf-rff,eqn:scalable-weights,eqn:streaming-weights}
are algebraically equivalent, we have arrived at
a streaming implementation of KAF with random features.

The general recommendation~\cref{eq:rfeatures} for the number $s$ of random features
may not be appropriate for the streaming setting because $s$ depends on the number
$n$ of training samples.  Our empirical work supports a more aggressive choice:
\begin{equation} \label{eqn:empirical-features}
s = \mathrm{Const} \cdot \ell.
\end{equation}
In other words, the number $s$ of features can be proportional to the
dimension $\ell$ of the regression model, which is chosen in advance.

\subsection{Streaming PCA}
\label{sec:nys}

To complete the description of our streaming KAF algorithm, we must provide an efficient method for computing
a low-rank approximation of the feature covariance matrix $\mtx{C}_{xx}$ appearing in~\cref{eqn:covariance-stream}.

Evidently, $\mtx{C}_{xx}$ is the covariance of vectors that are presented
to us sequentially.  Therefore, the low-rank approximation $\lowrank{\mtx{C}_{xx} + \mu\Id}_{\ell}$
amounts to a streaming PCA problem.  We will perform this computation using the
randomized Nystr{\"o}m method~\cite{2011Halko,Git13:Topics-Randomized,LLS+17:Algorithm-971,2017Tropp,martinsson2020randomized};
see~\cref{sec:related} for a short discussion of alternatives.

The Nystr{\"o}m approximation of a positive-semidefinite (psd) matrix $\mtx{C} \in \RR^{s \times s}$
with respect to a test matrix $\mtx{\Omega} \in \RR^{s \times k}$ is the best psd approximation
with the same range as $\mtx{C\Omega}$.  The construction dates back to the early literature
on integral equations~\cite{Nys30:Uber-Praktische}; it is intimately connected to Schur complements and Cholesky
factorization.  The randomized Nystr{\"o}m approximation involves a test matrix $\mtx{\Omega}$
chosen at random.

We can implement randomized Nystr{\"o}m approximation in the streaming setting~\cite{2017Tropp}.
Draw and fix a random matrix $\mtx{\Omega} \in \RR^{s \times 2\ell}$ from the standard normal distribution.%
\footnote{It is important that the random matrix
$\mtx{\Omega}$ has $2\ell$ columns, not merely $\ell$.}
Instead of forming $\mtx{C}_{xx}$ as in~\cref{eqn:covariance-stream},
we compute the product $\mtx{B} = \mtx{C}_{xx} \mtx{\Omega} \in \RR^{s \times \ell}$ via the iteration
$$
\mtx{B} = \mtx{0}_{s \times \ell}
\quad\text{and}\quad
\mtx{B} \gets \mtx{B} + \varphi(\mtx{x}_i) (\varphi(\mtx{x}_i)^\tsp \mtx{\Omega}).
$$
After we have streamed all of the data, we carefully%
\footnote{Do not use the formula~\cref{eqn:Cxx-Nystrom} as written!  See~\cref{alg:nys}.}
form a Nystr{\"o}m
approximation of the covariance 
and extract its eigenvalue decomposition:
\begin{equation} \label{eqn:Cxx-Nystrom}
\check{\mtx{C}}_{xx} \coloneqq \mtx{B} (\mtx{\Omega}^* \mtx{B})^\dagger \mtx{B}^*
	= \mtx{Q\Lambda Q}^\tsp.
\end{equation}
The randomized Nystr{\"o}m approximation $\check{\mtx{C}}_{xx}$
provides a good low-rank approximation of the covariance $\mtx{C}_{xx}$; see~\cite[Thms.~4.1--4.2]{2017Tropp}.
Our ultimate formula for the weight matrix becomes
\begin{equation} \label{eqn:fast-weights}
\check{\mtx{W}}_{q, \ell} = \mtx{C}_{gx} \cdot (\lowrank{ \check{\mtx{C}}_{xx} + \mu \Id }_{\ell})^{\dagger}.
\end{equation}
We can easily complete this computation because we have the eigenvalue decomposition
of the approximation $\check{\mtx{C}}_{xx}$ at hand.  The final target function becomes
$\check{f}_{q, \ell}(\mtx{y}) \coloneqq \check{\mtx{W}}_{q,\ell} \cdot \varphi(\mtx{y})$.

\cref{alg:nys} provides numerically stable pseudocode for the
randomized Nystr{\"o}m method applied to a sequence of random features.
This method is based on~\cite{LLS+17:Algorithm-971,2017Tropp}.

Using ordinary Gaussian random features, the arithmetic cost of forming the matrix
$\mtx{B}$ is $O((\ell + d)sn)$.
The algorithm uses auxiliary arithmetic $O(\ell^2 s)$,
and the storage requirement is just $O(\ell s)$.

\begin{remark}[Powering] \label{rem:powering}
The randomized Nystr{\"o}m method always underestimates the eigenvalues
of the covariance matrix.  If necessary, we can reduce this effect by incorporating
powering or Krylov subspace techniques~\cite{LLS+17:Algorithm-971,martinsson2020randomized}.
In the streaming setting, these modifications require us to construct and store
the full covariance matrix $\mtx{C}_{xx}$.  In our numerical work,
these refinements did not improve the quality of forecasting,
but they may merit further study.
\end{remark}

\begin{algorithm}[t] 
  \caption{\textsl{Randomized Nystr{\"o}m for featurized data~\cite[Sec.~19.4.3]{martinsson2020randomized}.}
  See \cref{sec:nys}.}
  \label{alg:nys}
  
  \vspace{0.5pc}
  
  Given a random feature map $\texttt{feat}$ and a data matrix $\mtx{X} \in \RR^{d \times n}$,
  this procedure computes an $\ell$-truncated eigenvalue decomposition $\mtx{Q \Lambda Q}^\tsp$
  of the covariance $\mtx{C}_{xx} = \varphi(\mtx{X}) \varphi(\mtx{X})^\tsp$ of the featurized data using the randomized
  Nystr{\"o}m method with $2\times$ oversampling. 
  
  \begin{algorithmic}[1]

	\vspace{1pc}

	\Function{FeatNystr{\"o}m}{RFF \texttt{feat}, $\mtx{X} \in \RR^{d \times n}$, $\ell \in \mathbb{N}$)} 

	\State	$\mtx{Q} \gets \texttt{orth}( \texttt{randn}(s, 2 \ell) )$
		\Comment	Random subspace, oversampling $\ell \to 2\ell$
		\State	$\mtx{Z} \gets \texttt{feat}.\textsc{MultCov}(\mtx{X}, \mtx{Q})$
			\Comment	Stream the product $\varphi(\mtx{X}) \varphi(\mtx{X})^\tsp \mtx{Q}$

	\State	$\nu \gets \texttt{eps}(\norm{\mtx{Z}}_{\mathrm{F}})$
		\Comment	Compute shift	
	\State	$\mtx{Z} \gets \mtx{Z} + \nu \mtx{Q}$
		\Comment	Shift for stability
	\State	$\mtx{T} \gets \texttt{chol}(\mtx{Q}^\tsp \mtx{Z})$
		\Comment	Upper-triangular Cholesky factorization
	\State	$\mtx{S} \gets \mtx{Z} / \mtx{T}$
		\Comment	Solve triangular systems

	\State	$(\mtx{Q}, \mtx{\Sigma}, \sim) \gets \texttt{svd}( \mtx{S})$
		\Comment	Compact SVD
	
	\State	$\mtx{\Lambda} \gets \max\{ \mathbf{0}, \mtx{\Sigma}^2 - \nu \Id \}$
		\Comment 	Remove shift to get eigenvalues
	
	\State	$\mtx{Q} \gets \mtx{Q}(:, 1:\ell)$ and $\mtx{\Lambda} \gets \mtx{\Lambda}(1:\ell, 1:\ell)$
		\Comment	Truncate to rank $\ell$
					
	\State	\Return		$(\mtx{Q} \in \RR^{s \times \ell}, \mtx{\Lambda} \in \RR^{\ell \times \ell})$
	\EndFunction

\vspace{0.5pc}

\end{algorithmic}
\end{algorithm}

\subsection{Other observables} \label{sec:gen-scalable-kaf}

We can easily extend streaming KAF to the more general setting outlined
in~\cref{sec:other-observables}.  Suppose we wish to use a general
covariate $u : \mathcal{M} \to \RR^{d'}$ to predict a general response
variable $g : \mathcal{M} \to \RR^r$ after $q$ time steps.
Let $\tilde{\kappa} : \RR^{d' \times d'} \to \RR_+$ be a positive-definite kernel
on the covariate space, with associated feature map $\tilde{\varphi} : \RR^{d'} \to \RR$.

To train, we acquire data in the form of measured values of the covariate paired with measured values of the
lagged response: $(\mtx{u}_i, \mtx{g}_{q+i})$
where $\mtx{u}_i = u(\mtx{x}_i)$ and $\mtx{g}_i = g(\mtx{x}_i)$ for $i = 0, \dots, n-1$.
In this setting, the underlying state trajectory $(\mtx{x}_0, \dots, \mtx{x}_{n-1})$ is unknown.
By streaming the observable data, we compute the matrices
$$
\mtx{C}_{uu} \gets \mtx{C}_{uu} + \varphi(\mtx{u}_i) \varphi(\mtx{u}_i)^\tsp
\quad\text{and}\quad
\mtx{C}_{gu} \gets \mtx{C}_{gu} + \mtx{g}_{i+q} \varphi(\mtx{u}_i)^\tsp.
$$
Finally, we determine the weights:
$$
\mtx{W}_{q, \ell} \coloneqq \mtx{C}_{gu} \cdot (\lowrank{\mtx{C}_{uu} + \mu \Id}_{\ell})^\dagger.
$$
As before, the randomized Nystr{\"o}m method serves for the streaming PCA computation.

Now, suppose that we observe a covariate $\mtx{v} \in \RR^{d'}$, where $\mtx{v} = u(\mtx{y})$ for an unknown state $\mtx{y} \in \mathcal{M}$.
We forecast the lagged response $g(F^q(\mtx{y}))$ as 
\begin{equation} \label{eqn:gen-forecast}
\hat{g}_{q, \ell}(\mtx{v}) \coloneqq \mtx{W}_{q, \ell} \cdot \tilde{\varphi}(\mtx{v}).
\end{equation}
Our approach gives a principled approximation of the optimal forecast of the response
given the observed covariate, as described in~\cref{sec:koop}.

\subsection{Resource usage}
\label{sec:rff-kaf-cost}

\Cref{alg:rff-kaf} lists pseudocode for the general streaming KAF method  outlined in~\cref{sec:gen-scalable-kaf}.
\Cref{tab:kaf} compares the costs against a na\"ive implementation of KAF.  We also list the costs of
streaming KAF with fast random features (Fast Streaming KAF; see~\cref{rem:fast-features}), omitting an exposition.

First, we discuss the costs of the training step of streaming KAF
with covariate data $\mtx{X} \in \RR^{d' \times n}$ and (lagged) observable data $\mtx{G} \in \RR^{r \times n}$.
Assume that the truncation rank $\ell \leq s$, where $s$ is the number of random features.

To construct random feature descriptors, we draw and store $O(d's)$ normal random variables.
The Nystr{\"o}m approximation of the featurized covariance matrix
involves $O((d'+\ell) s n)$ arithmetic and local storage $O(\ell s)$.
The covariate--response matrix requires $O((d' + r)sn)$ arithmetic
and storage $O(rs)$.
To form the prediction weights, we expend $O(\ell r s)$ arithmetic
and $O(rs)$ storage.
In practice, the Nystr{\"o}m approximation is the most expensive step.

The total storage required for the forecasting model consists of the $O(d's)$ storage
for the random feature descriptors and the $O(rs)$ storage for the
prediction weights.

In the forecasting step, we simply featurize the test data and
form a matrix--matrix product.  This step uses $O((d'+r)s)$ arithmetic
per initial condition (IC), but no additional storage.

Let us summarize.  In comparison with na{\"i}ve KAF, the streaming KAF method
is significantly faster because it is a streaming method.  The precise improvements
to storage and arithmetic costs depend on several parameters.  Loosely, the
streaming method reduces training arithmetic by a factor of about $n/s$ and
reduces training storage by a factor of about $n^2/s$.
For forecasting, the arithmetic and storage both decrease by a factor of $n/s$.

\begin{remark}[Implementation]
For reasons of modularity, the pseudocode and our prototype implementation take two passes over the data,
but they are mathematically equivalent to the streaming KAF algorithm.
\end{remark}

\begin{algorithm}[t!] 
  \caption{\textsl{Scalable Kernel Analog Forecasting.}
  Implements~\cref{sec:gen-scalable-kaf}.}
  \label{alg:rff-kaf}
  
  \vspace{0.5pc}

  The method \textsc{Train} takes covariate data $\mtx{X} \in \RR^{d \times n}$
  and (lagged) response data $\mtx{G} \in \RR^{r \times n}$ as input.
  It constructs a random feature map with parameters $(\gamma, d; s)$
  and builds a forecasting model for the response data $\mtx{G}$ using truncation rank
  $\ell \in \mathbb{N}$.  The method \textsc{Forecast} uses the model to make estimates
  of the response from the covariates listed as columns of $\mtx{Y} \in \RR^{d \times m}$.  
  
  \begin{algorithmic}[1]

	\vspace{1pc}
	\State	\textbf{local variables} RFF \texttt{feat}
		\Comment	Random feature map $\varphi:\RR^{d} \to \RR^s$
	\State	\textbf{local variables} $\mtx{W} \in \RR^{r \times s}$
		\Comment	Prediction weights
	
	\vspace{0.5pc}

	\Function{Train}{$\mtx{X} \in \RR^{d' \times n}$, $\mtx{G} \in \RR^{r \times n}$}

	\State	$\texttt{feat} \gets \textsc{RFF}(\gamma, d; s)$
		\Comment	Initialize random feature map
	\State	$(\mtx{Q}, \mtx{\Lambda}) \gets \textsc{FeatNystr{\"o}m}(\texttt{feat}, \mtx{X}; \ell)$
		\Comment	Factor $\lowrank{\varphi(\mtx{X})\varphi(\mtx{X})^\tsp}_{\ell}$; see \cref{sec:nys}
	\State	$\mtx{\Lambda} \gets \mtx{\Lambda} + \mu \max(\mtx{\Lambda}) \cdot \Id$
		\Comment	Filter eigenvalues; $\mu = 10^{-6}$
	\State	$\mtx{C} \gets \texttt{feat}.\textsc{RMultAdj}(\mtx{X}, \mtx{G})$
		\Comment	Form product $\mtx{G} \varphi(\mtx{X})^\tsp \in \RR^{r \times s}$
	\State	$\mtx{W} \gets ((\mtx{C}\mtx{Q}) / \mtx{\Lambda}) \mtx{Q}^\tsp$
		\Comment 	Compute prediction weights
	\EndFunction
	
	\vspace{0.5pc}
	
	\Function{Forecast}{$\mtx{Y} \in \RR^{d \times m}$}
	
	\State	$\hat{\mtx{F}} \gets \texttt{feat}.\textsc{RMult}(\mtx{Y}, \mtx{W})$
		\Comment	Form $\mtx{W} \,\varphi(\mtx{Y})$
	\State	\Return	$\hat{\mtx{F}} \in \RR^{d\times m}$
		\Comment Forecasts for columns of $\mtx{Y}$
	\EndFunction
		
\vspace{0.5pc}

\end{algorithmic}
\end{algorithm}

 \begin{table}[t]
    \centering
    \caption{\textbf{Resource usage for training and for a single forecast:} Covariate dimension $d'$, response dimension $r$, with $n$ training samples, 
    $s$ random features, truncation rank $\ell$.  Assumes $\ell \leq s \leq n$.  Constants are suppressed.  The fast streaming method uses a more efficient random feature construction.  See \cref{sec:kaf-cost,sec:rff-kaf-cost}.}
    \label{tab:kaf}
    \begin{tabular}{|c|c|c|c|c|}
        \hline
		\multicolumn{2}{|c|}{\,} & \textbf{Na{\"i}ve KAF} & \textbf{Streaming} & \textbf{Fast Streaming} \\
        \hline
		\multirow{3}{*}{\textbf{Training}} 
			& Streaming & $\times$ & $\checkmark$ & $\checkmark$ \\
			& Arithmetic & $d'n^2$ & $(\ell + d') s n + r \ell s$ & $(\ell + \log d')sn + r\ell s$  \\
			& Local storage & $n^2$ & $(\ell + r + d')s$ & $(\ell + r) s$ \\			
		\hline
		\multirow{2}{*}{\textbf{Forecast}}
			& Storage for model & $(r+d')n$ & $(d' + r)s$ & $rs$ \\
			& Arithmetic (per IC) & $rn$ & $(r+d')s$ & $(r + \log d')s$ \\
		\hline
    \end{tabular}
\end{table}

\section{Experiments}
\label{sec: experiments}

This section showcases experiments that demonstrate the practical performance of streaming KAF.
We study forecasting skill for several benchmark dynamical systems, we investigate
sensitivity to algorithm parameters, and we make comparisons with the na{\"i}ve implementation of KAF.
The code for reproducing the experiments is available as a supplement to this paper.

\subsection{The Lorenz models}
\label{sec:lorenz}

Our experiments focus on the Lorenz '63 model, a classical three-dimensional dynamical system known to exhibit chaotic behavior.
We also test the method on the two-phase Lorenz '96 model, a higher-dimensional system that has periodic, quasi-periodic,
and chaotic regimes.  This subsection summarizes the models and the parameters that give rise to different types of dynamics.

\subsubsection{Lorenz '63}

The Lorenz '63 (L63) model was introduced by Edward Lorenz in 1963 as a crude model of atmospheric convection~\cite{Lorenz63}.
Although this example is simple, its properties have been studied extensively, and it is known to exhibit many of the features that make forecasting challenging in more complex systems, including fractal attractors~\cite{Tucker99} and mixing dynamics~\cite{LuzzattoEtAl05}. 

The L63 model is defined via the following system of differential equations.
For a state $\mtx{x} = (x_1, x_2, x_3) \in \RR^3$,
\begin{equation} \label{eqn:L63}
\begin{aligned}
\dot{\mtx{x}}(t) &= \mtx{V}(\mtx{x}(t)) \quad\text{with initial condition $\mtx{x}(0) = \mtx{x}_{\mathrm{init}}$}; \\ 
 V_1(\mtx{x}) &= \sigma(x_2 - x_1); \quad
 V_2(\mtx{x}) = x_1(\mu - x_3); \quad
 V_3(\mtx{x}) = x_1x_2 - \beta x_3.
\end{aligned}
\end{equation}
The classical parameters for the L63 system that generate chaotic dynamics are $(\sigma,\mu,\beta) = (10,28,8/3)$.
This choice leads to the famous ``butterfly attractor,'' a compact set in $\RR^3$ 
with fractal dimension $\approx 2.06$ that supports an ergodic invariant measure with Lyapunov
exponent $\lambda \approx 0.91$; see~\cite{Sprott03}.  \Cref{fig:lorenz} presents an illustration.

\subsubsection{Lorenz '96}

\begin{figure}[t]
    \centering
    \includegraphics[width=.45\textwidth,height=2.5in]{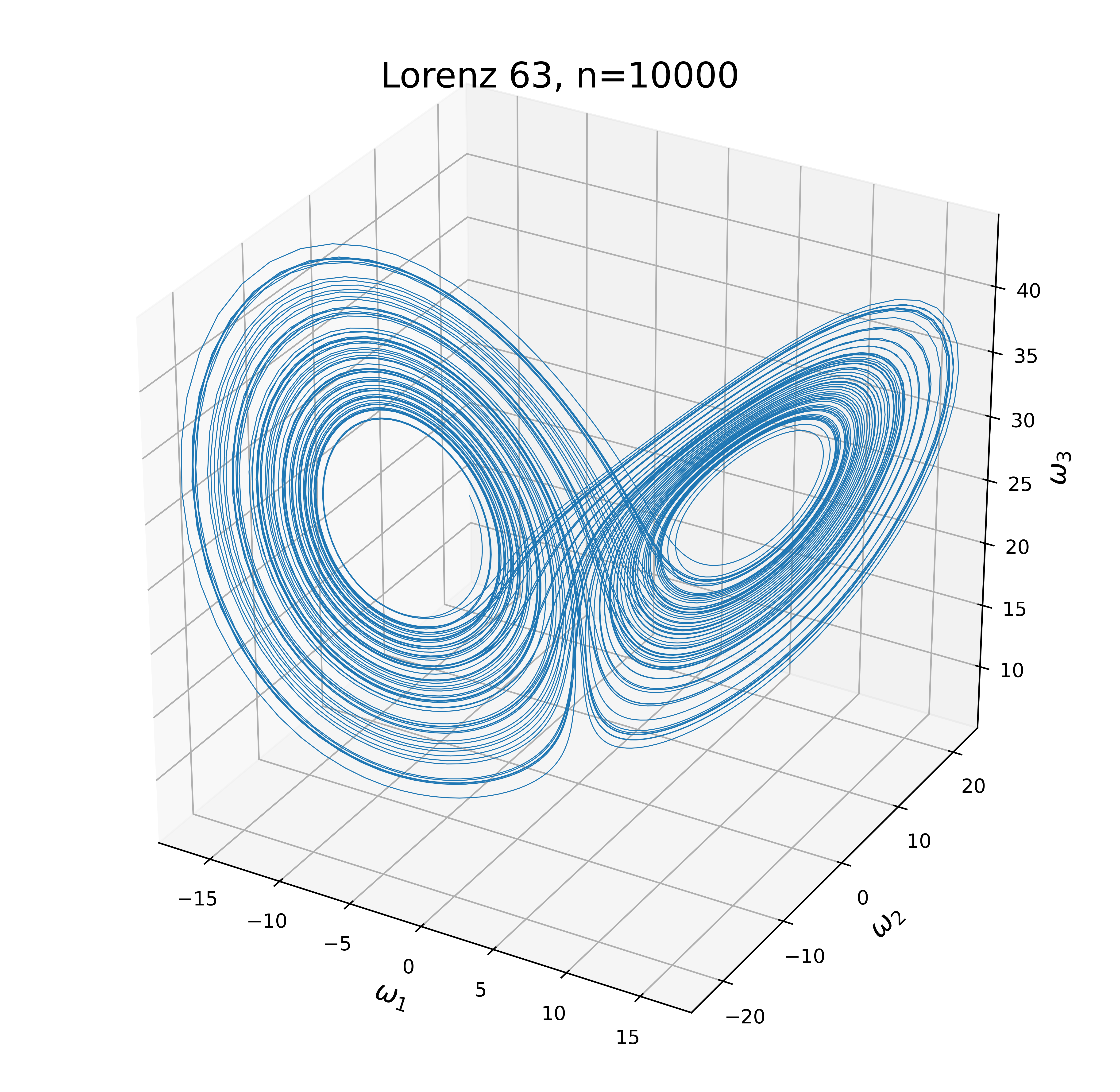} \hspace{1pc}
    \includegraphics[width=.45\textwidth]{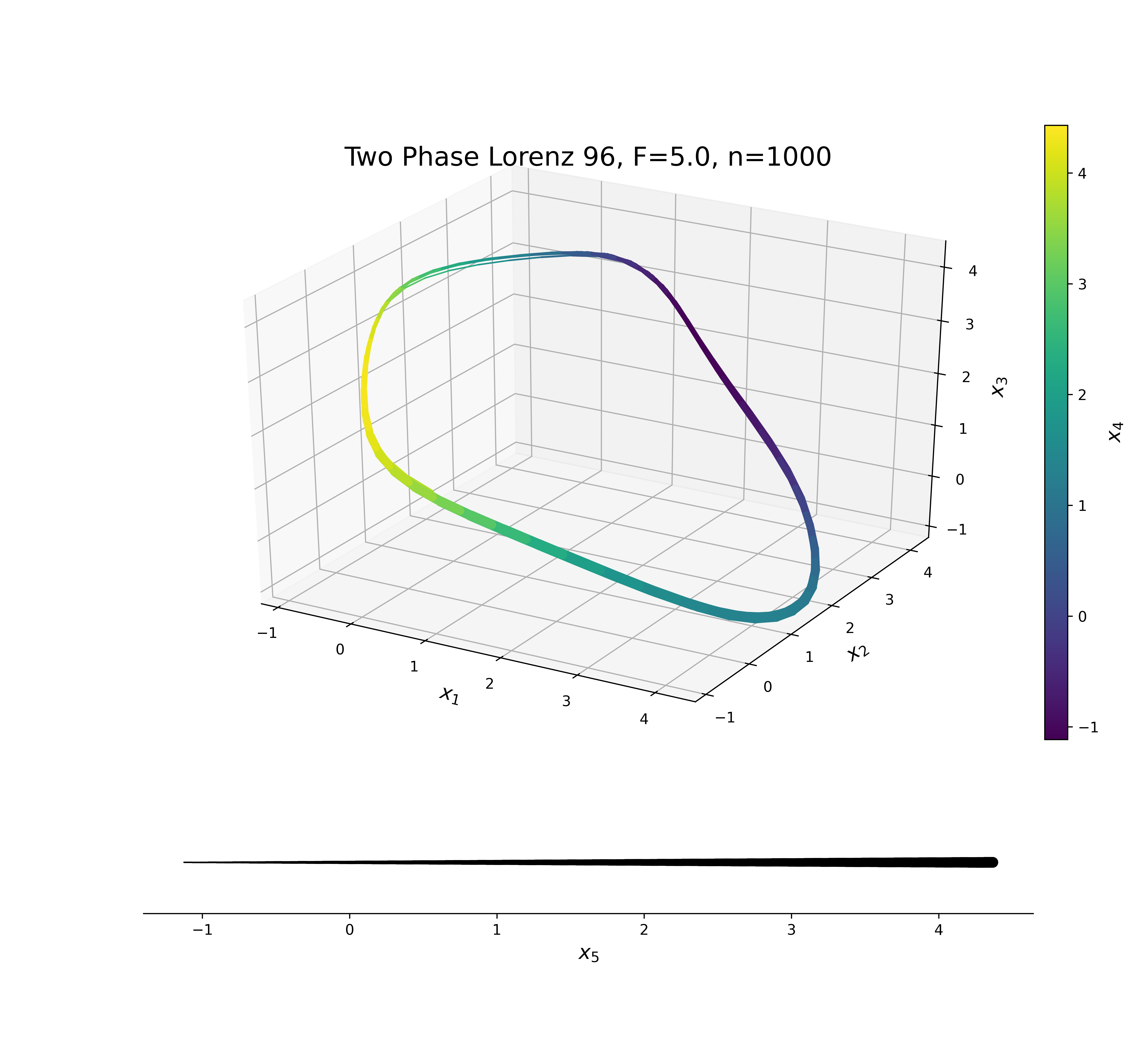} \\[1pc]
    \includegraphics[width=.45\textwidth]{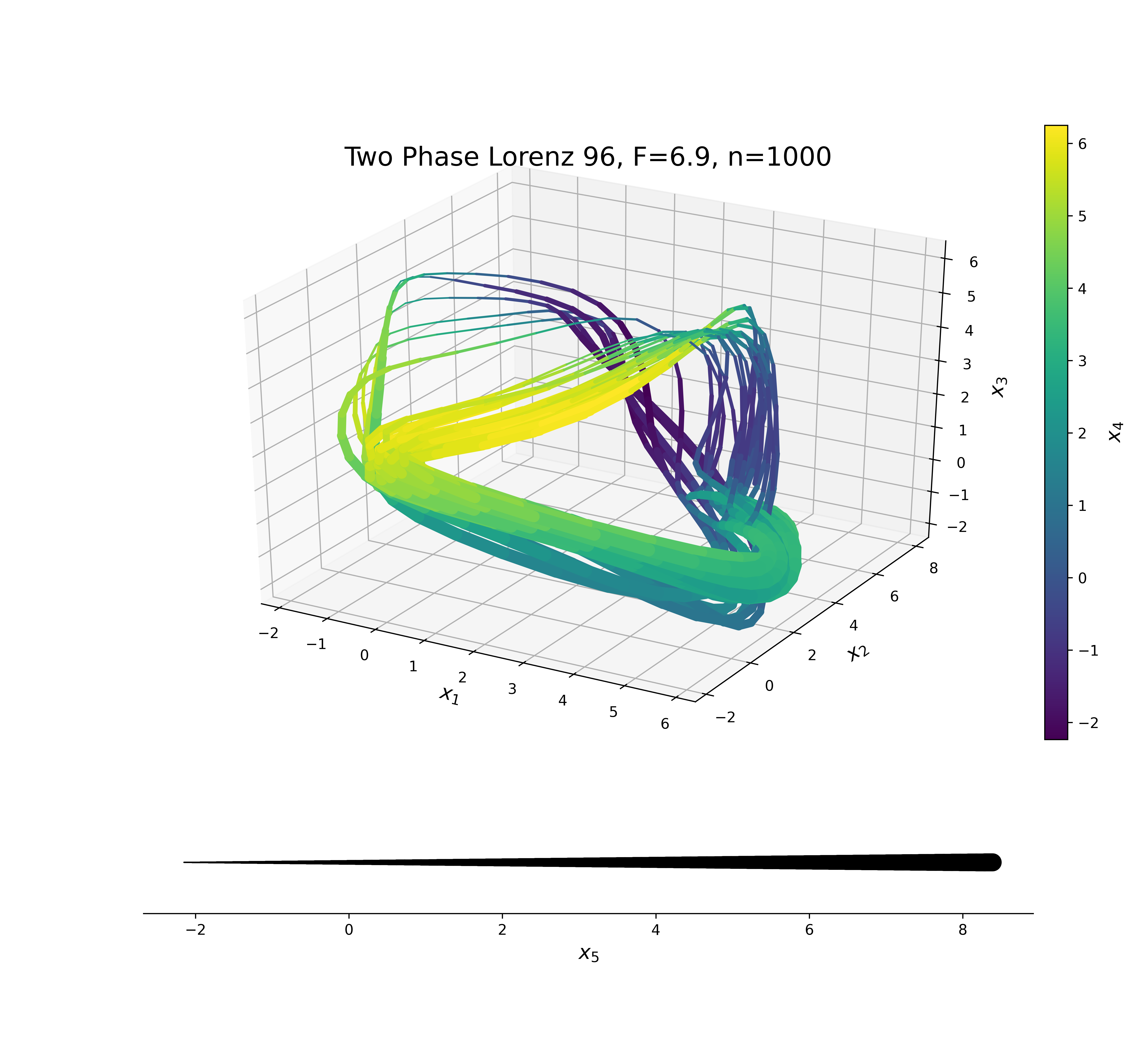} \hspace{1pc}
    \includegraphics[width=.45\textwidth]{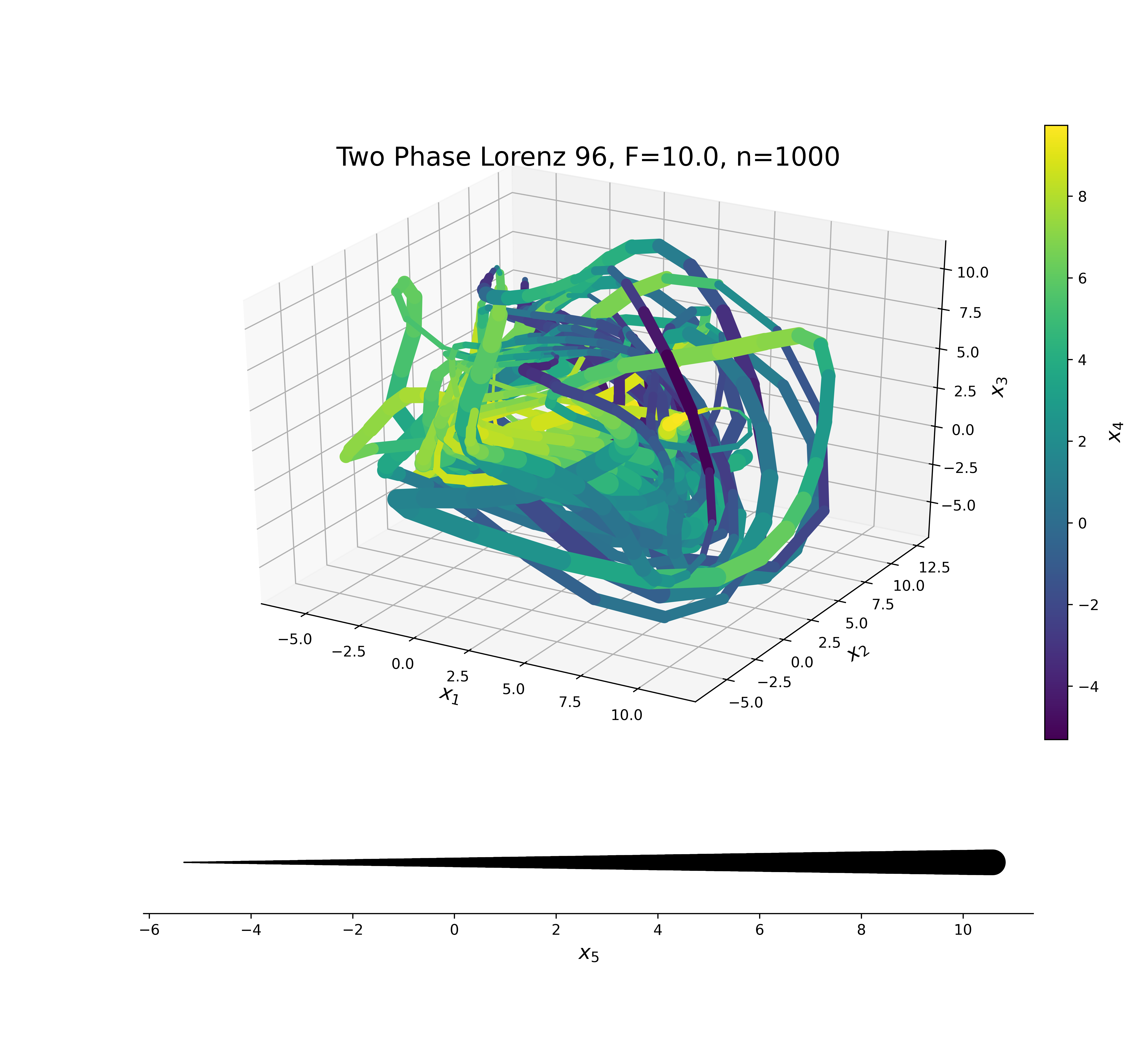}
    \caption{\textbf{Lorenz models.}  [top left] The L63 system in the chaotic regime.  [Other panels] Five slow dimensions of the L96 system.  The fourth dimension is plotted in color, and the fifth dimension is plotted as linewidth.  [top right] Periodic regime ($F = 5$).  [bottom left] Quasi-periodic regime ($F = 6.9$).  [bottom right] Chaotic regime ($F = 10$).  See~\cref{sec:lorenz} for details.}
    \label{fig:lorenz}

\end{figure}

We also consider the two-phase Lorenz '96 system (L96), as introduced in  \cite{1996Lorenz,2004Fatkullin}. 
This model has dynamics that occur on two distinct timescales, a set of ``slow variables" $\mtx x = \{x(k)\}_{k \in [K]}$ and a set of ``fast variables'' $\mtx z = \{z(j,k)\}_{j\in[J], k\in[K]}$. 
These variables evolve according to the following system of equations~\cite{2004Fatkullin}.  The boundary conditions $x(k+K) = x_k$ and $z(j,k+K) = z(j,k)$ for $k \in [K]$ and $z(j+J, k) = z(j, k+1)$ for $j \in [J]$; the dynamics are
\begin{equation} \label{eqn:L96}
\begin{aligned}
\dot{x}(k) &= -x(k-1) \left( x(k-2) - x(k+1) \right) - x(k) + F + \frac{h_x}{J} \sum\nolimits_{j=1}^J z(j,k)\\
\dot{z}(j,k) &= \frac{1}{\varepsilon}\left( -z(j+1,k) \left( z(j+2,k) - z(j-1,k) \right) - z(j,k) + h_y \cdot x(k) \right).
\end{aligned}
\end{equation}
As in~\cite{2004Fatkullin}, we set the parameters $(h_x, h_y, K, J, \varepsilon) = (-0.8, 1, 9, 8, 1/128)$.
Depending on the value of the forcing constant $F$, three distinct regimes of behavior emerge.
\begin{itemize}
\item $F = 5$ yields a periodic system;
\item $F = 6.9$ yields a quasi-periodic system; and
\item $F = 10$ yields a fully {chaotic} system.
\end{itemize}
See~\cref{fig:lorenz} for typical trajectories. 

In our experiments, we seek to forecast the future values of the slow variables $x(1), \dots, x(9)$
of the coupled system \textit{using only the slow variables as input data}.  This setup is motivated by the experiments of \cite{BurovEtAl21}, which studied KAF for multi-scale systems but did not investigate the scalability as a function of the amount of training data.

\subsection{Experimental setup}

All of our experiments are performed using data obtained by integrating the governing equations of the L63 and L96 systems.
Here are the details about how we apply streaming KAF to make forecasts and evaluate the results.

\begin{itemize}
\item For L63, the training data consists of states $\mtx x_1, \mtx x_2, \dots, \mtx x_{n} \in \mathbb{R}^d$ generated by iterating the dynamics:
$$
\mtx x_j \coloneqq F(\mtx x_{j-1}), \quad \quad j = 1, 2, \dots, n-1,
$$
where $F$ is the flow map obtained by discretizing~\cref{eqn:L63} with time step $\mathrm{d}t = .01$.   
 
\item For L96, we first generate the full $81$-dimensional system of slow  and fast variables:
$$
[\mtx x_j, \mtx z_j]  \coloneqq F([\mtx x_{j-1}, \mtx z_{j-1}]), \quad \quad j = 1, 2, \dots, n-1,
$$
where $F$ is the flow map obtained by discretizing~\cref{eqn:L96} with time step ${\diff}t = .01$. 
We then  form the training matrix $\mtx X = [\mtx x_1, \dots, \mtx x_n] \in \mathbb{R}^{9 \times n}$ using only the slow variables. 

\item To evaluate the performance, we use the normalized root mean square error (RMSE) metric for the forecast error.  For a single response variable $i^*$, consider the test set $\mtx{Y}_{i^*}$ and the true trajectory $\mtx{Y}_{q, i^*}$:
$$
\begin{aligned}
\mtx Y_{i^{*}} &= [ \mtx y_0(i^{*}), \mtx y_1(i^{*}), \dots, \mtx y_{m-1}(i^{*}) ] \in \mathbb{R}^{1 \times m} \\
\mtx Y_{q, i^{*}} &= [ \mtx y_q(i^{*}), \mtx y_{q+1}(i^{*}), \dots, \mtx y_{q+m-1}(i^{*}) ] \in \mathbb{R}^{1 \times m}.
\end{aligned}
$$
For the forecast $f_{q, \ell, i^*}$ of the response variable, applied columnwise, we define the error

\begin{equation}
\label{eq:rmse}
\text{RMSE}(f_{q, \ell, i^*}(\mtx Y)) = \frac{\| f_{\ell, q, i^*} - \mtx{Y}_{q, i^{*}} \|_2}{\sqrt{m} \cdot \operatorname{std}(\mtx Y_{q, i^{*}})} 
\end{equation}
where $\operatorname{std}({\mtx z})$ denotes the standard deviation of the vector $\mtx z$.

\item	For each system and each set of parameter specifications, we consider $5$ sets of tests $\mtx Y_1, \dots, \mtx Y_5$, each with $m = 10,000$ data points (columns) of the same form as the training data.  The first test data set $\mtx{Y}_1$ is obtained by evolving the system from the final point $\mtx{x}_n$ in the training data.  For the remaining test sets, the initial condition is the final point in the previous set.  In all figures, the line series represents the average of the errors resulting from each of the $5$ tests, and the shaded region around the error lines represents one standard deviation of uncertainty around the average.

\item In each of the plots presented in sections \ref{sec:l63} and \ref{sec:l96}, the kernel inverse bandwidth $\gamma$ and dimension of regression model $\ell$ are fixed as the size of the training data $n$ increases.  For any particular plot in these sections, the values of $\gamma$ and $\ell$ were chosen based on a minimal amount of manual tuning at fixed training sample size $n=10,000$.  As such, the corresponding error curves level off as $n$ is increased from $10,000$ to $50,000$. 
Principled approaches to setting the inverse bandwidth parameter $\gamma$ and dimension of regression model $\ell$ are discussed in sections \ref{sec:bandwidth-param} and \ref{sec:dim-param}, respectively. 

\item Table  \ref{tab:L63_train_time} illustrates that gently increasing the inverse bandwidth $\gamma$ and regression model dimension $\ell$ together as the size of the training data $n$ increases serves as a good rule of thumb for improving the streaming KAF accuracy with increasing $n$.  Table \ref{tab:small_s} indicates that the number of random features $s$ can be taken to be proportional to $\ell$, resulting in faster forecasting and incurring only a small loss in accuracy.

\item	As presented in \cref{alg:rff-kaf}, streaming KAF is implemented with two passes over the training data, but it is algebraically equivalent to a true streaming method.  We use ordinary Gaussian random features (rather than the ``fast'' variant).  All the loops in~\cref{alg:rff} are vectorized with blocks of $1,000$ vectors.  We employ the randomized Nystr{\"o}m method described in~\cref{alg:nys}.

\item	The algorithms were implemented using the MATLAB programming language.  All data was collected on a MacBook Pro with 16 GB of RAM and with an 8-Core Intel Core i9 Processor, clocked at 2.3 GHz.
\end{itemize}

\subsection{Interpreting the results}

The normalized RMSE~\cref{eq:rmse} provides a measure of the quality of the forecast.  When the normalized RMSE reaches $1$, the expected square error is equal to the standard deviation of the response observable with respect to the invariant measure, and
the forecast is no longer providing useful information.  

In dynamical systems, the maximal Lyapunov exponent of a system is commonly used to summarize the level of ``unpredictability.'' 
The paper \cite{2019Wang} describes the intuitive meaning of this exponent:
``For a chaotic trajectory, an infinitesimal perturbation in the evolution gives rise to exponential divergence---the Lyapunov exponent expresses the rate of divergence.''
Hence, Lyapunov time is frequently used as a \textit{horizon for forecasting}.  Typically, a forecast is classified as ``good'' if the normalized RMSE only approaches 1 after several Lyapunov timescales.

For the L63 system, the Lyapunov exponent $\lambda \approx 0.91$ \cite{Sprott03}.  Thus, we can expect to make nontrivial forecasts of the state vector for several time units. 
On all L63 plots, we mark the Lyapunov time scale as a yardstick.  Note that there are other observables that remain predictable for much longer than the coordinates of the state vector~\cite{Giannakis21a}.

\subsection{Case study: L63}
\label{sec:l63}

\begin{figure}[t]
    \centering
      \includegraphics[width=.45\textwidth]{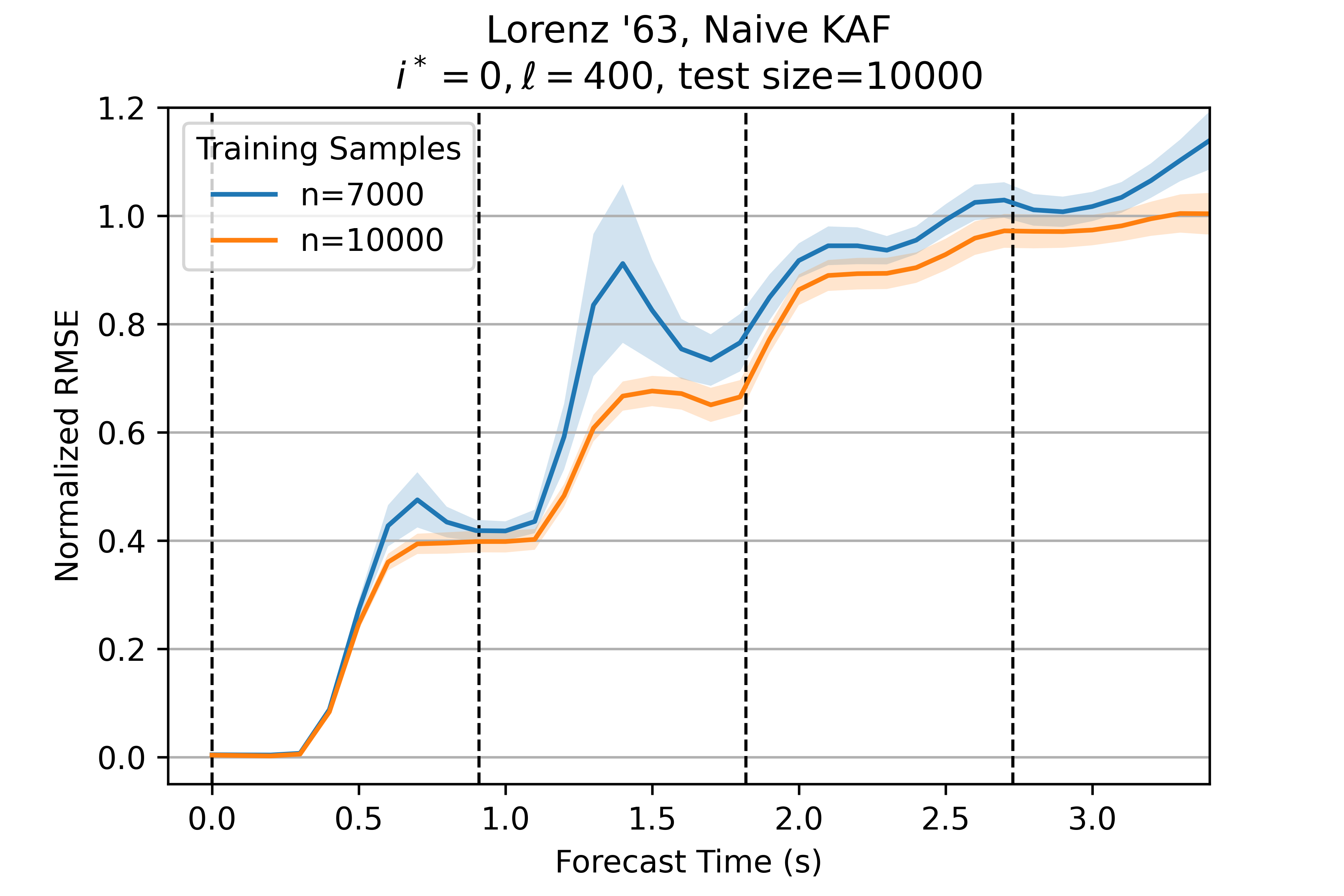} \quad
  \includegraphics[width=.45\textwidth]{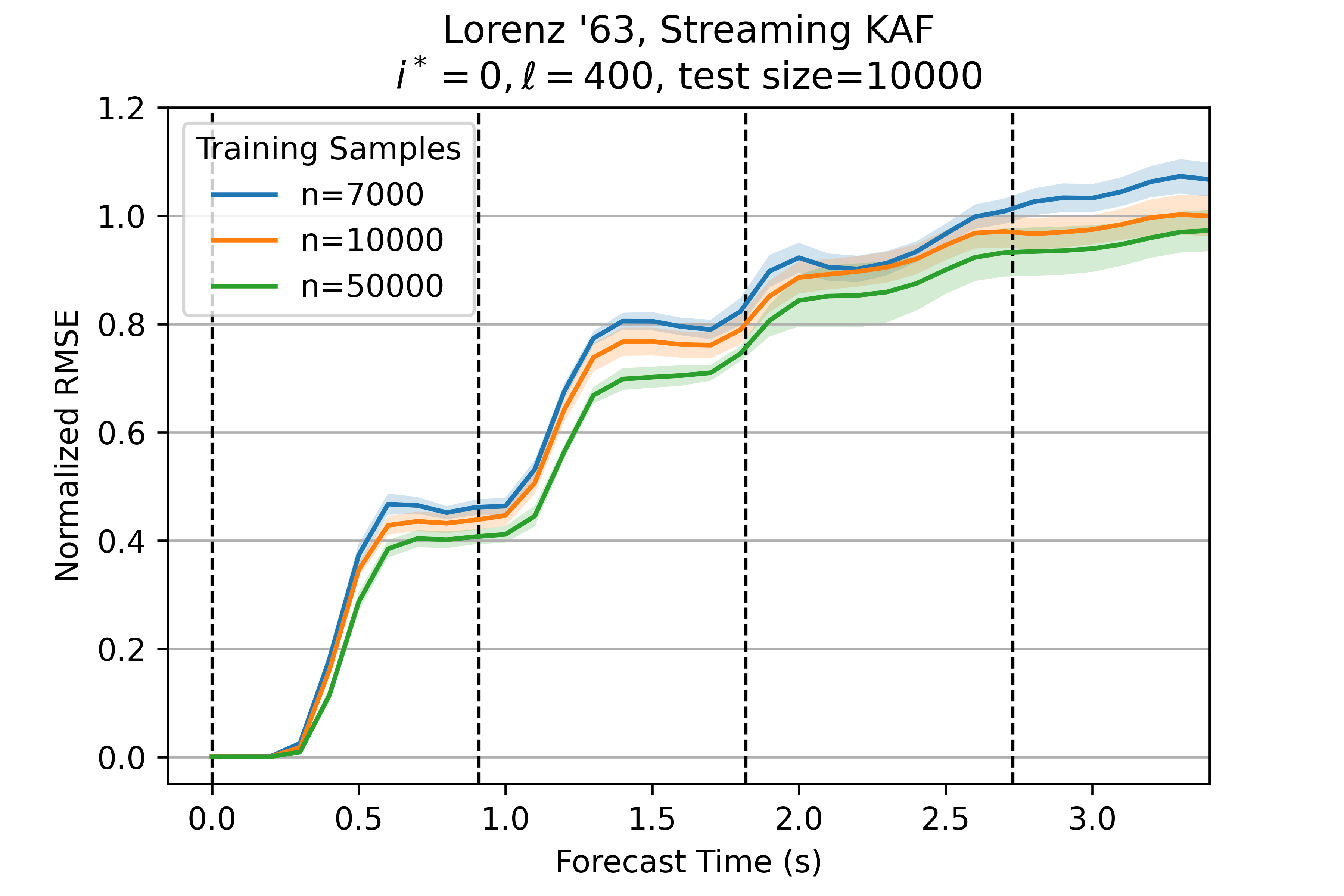} 
    \caption{\textbf{Lorenz '63: Forecast error versus amount of training data.}  Average normalized RMSE for forecasting the first state variable of L63 via na{\"i}ve KAF [left] and streaming KAF [right] as a function of the number $n$ of training points.
    The regression model has dimension $\ell = 400$, the kernel inverse bandwidth $\gamma = .05$, and the number of features $s = \sqrt{n} \log(n)$.  For $n=50,000$, na{\"i}ve KAF fails because of its computational cost. See~\cref{sec:l63}.}
    \label{fig:l63-rmse-vs-n}
\end{figure}

Our first experiment compares the forecasting skill of na{\"i}ve KAF and scalable KAF for the L63 system.
\Cref{fig:l63-rmse-vs-n}
explores how forecasts of the first state coordinate $i^{*} = 1$ improve as the number $n$ of training samples increases.
With $n=10,000$, both methods provide good predictions, with streaming KAF slightly better than na{\"i}ve KAF.  In particular, both algorithms can make informative forecasts over several Lyapunov time intervals.  As we will discuss in \cref{sec:timings}, the streaming method is far more efficient, and the na{\"i}ve method was unable to construct a forecasting model when $n = 50,000$.

\begin{figure}[t]
    \centering
  \includegraphics[width=.45\textwidth]{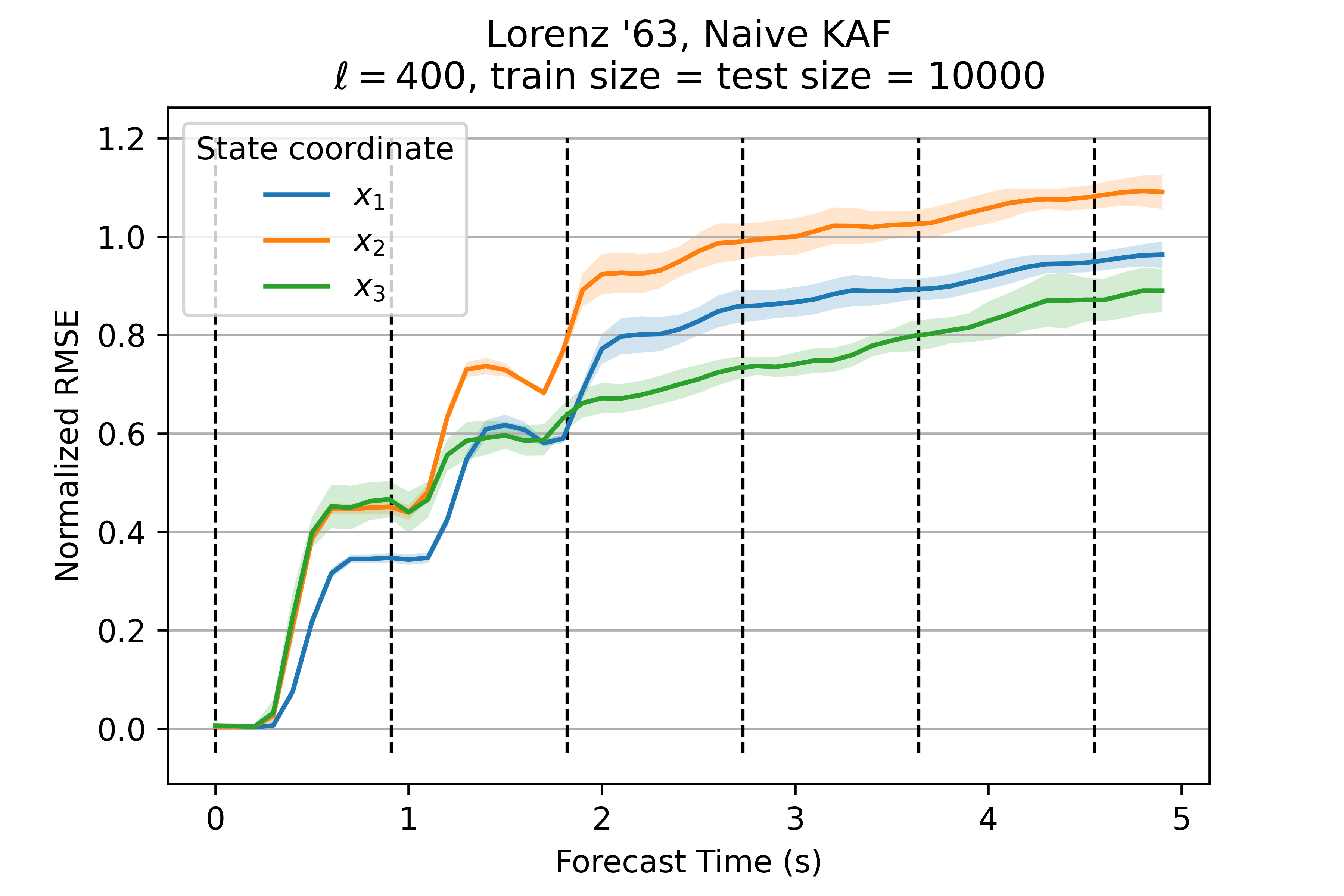} \quad 
  \includegraphics[width=.45\textwidth]{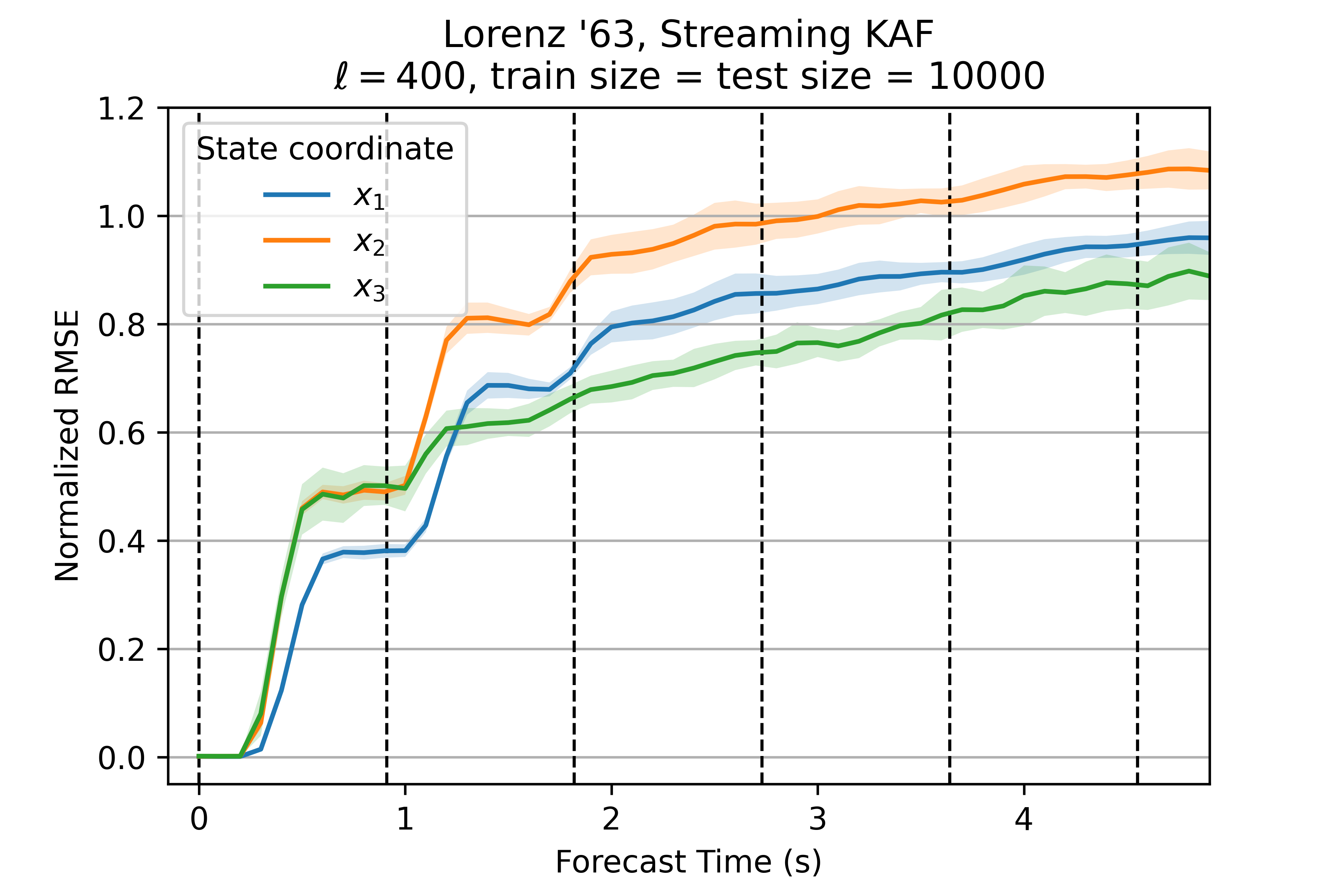}
    \caption{\textbf{Lorenz '63: Forecasting all three state variables.}  Average normalized RMSE for forecasting all three state coordinates $(x_1,x_2,x_3) = (\text{blue, orange, green})$ via na{\"i}ve KAF [left] and streaming KAF [right] with $n = 10,000$ training points.
    The regression model has dimension $\ell = 400$, the kernel inverse bandwidth $\gamma = .05$, and the number of features $s = \sqrt{n} \log(n)$.  See~\cref{sec:l63}.}
    \label{fig:l63-rmse-all}
\end{figure}

The KAF methodology has similar success at forecasting all three state variables.
For each of the three variables and with $n = 10,000$ training samples,
\cref{fig:l63-rmse-all} compares the forecasting error attained by
the na{\"i}ve and streaming methods.

\subsection{Case Study: L96}
\label{sec:l96}

In our second set of experiments,
we explore the performance of scalable KAF for the L96 system
in the periodic, quasi-periodic, and chaotic regimes documented in~\cite{BurovEtAl21}.
An increase in the forcing constant $F$ generates more chaotic behavior and,
unsurprisingly, reduces the time horizon for which KAF can make informative forecasts.

\begin{figure}[t]
    \centering
\includegraphics[width=.45\textwidth]{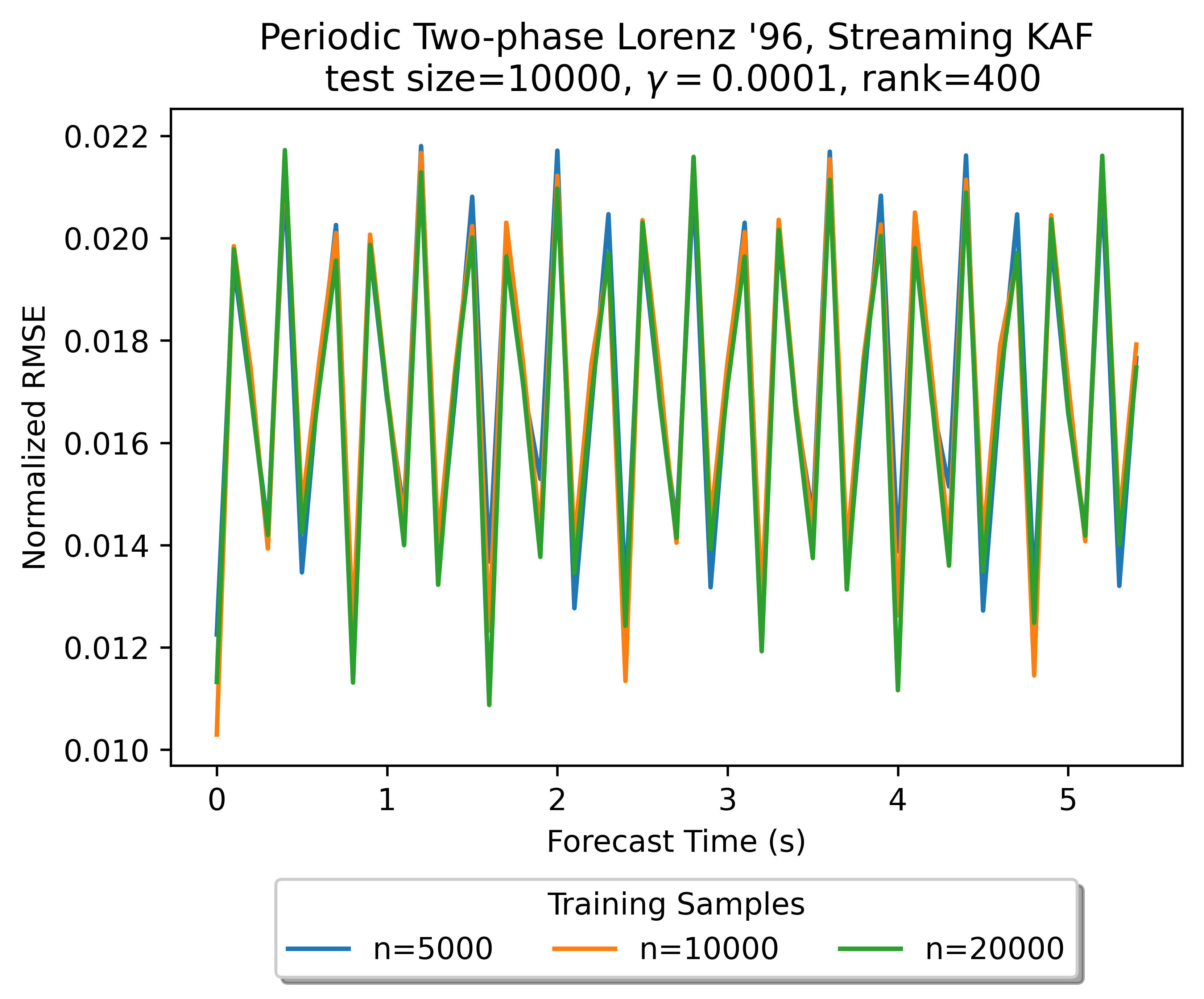} \\[1pc]
\includegraphics[width=.45\textwidth]{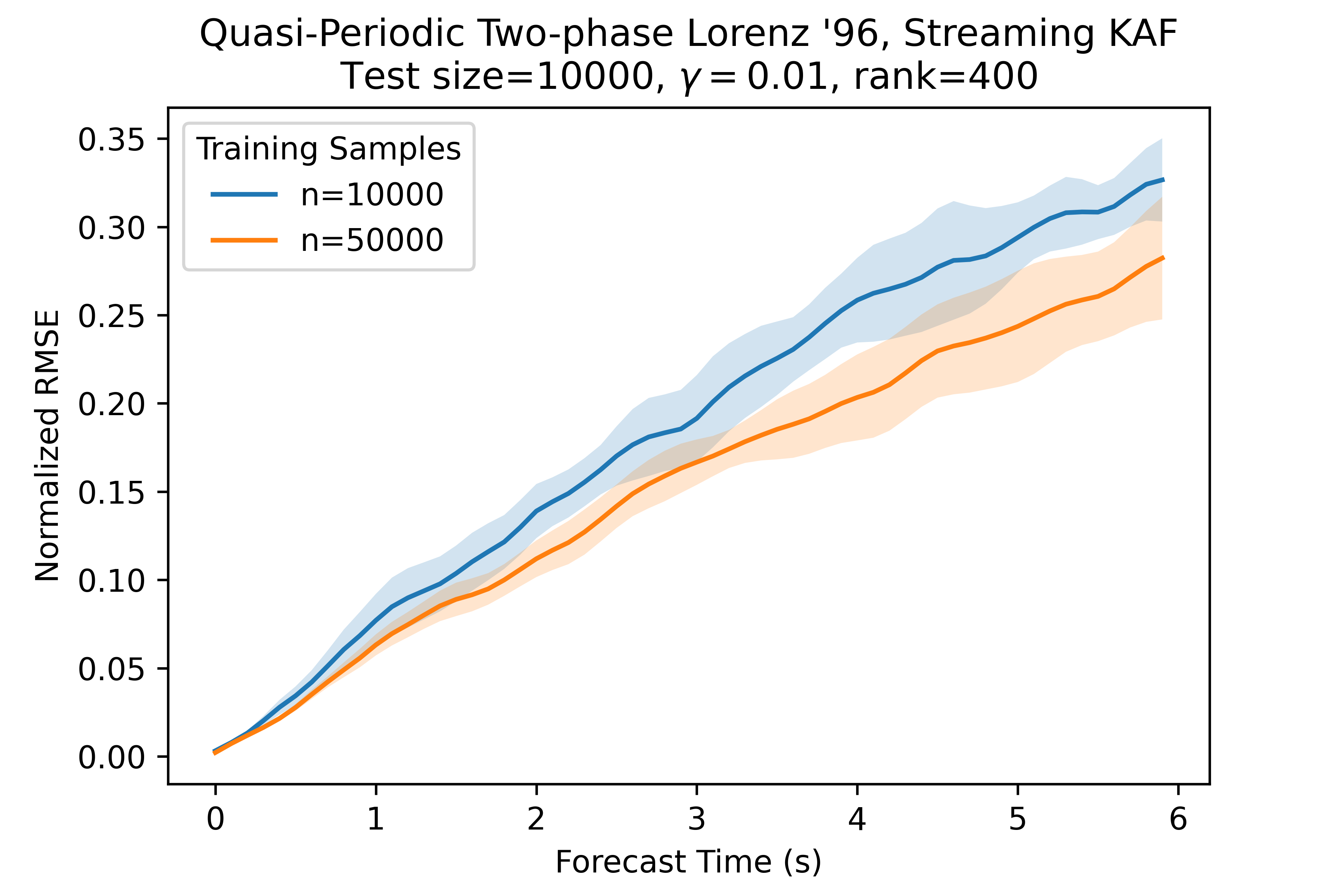}  \hspace{1pc}
\includegraphics[width=.45\textwidth]{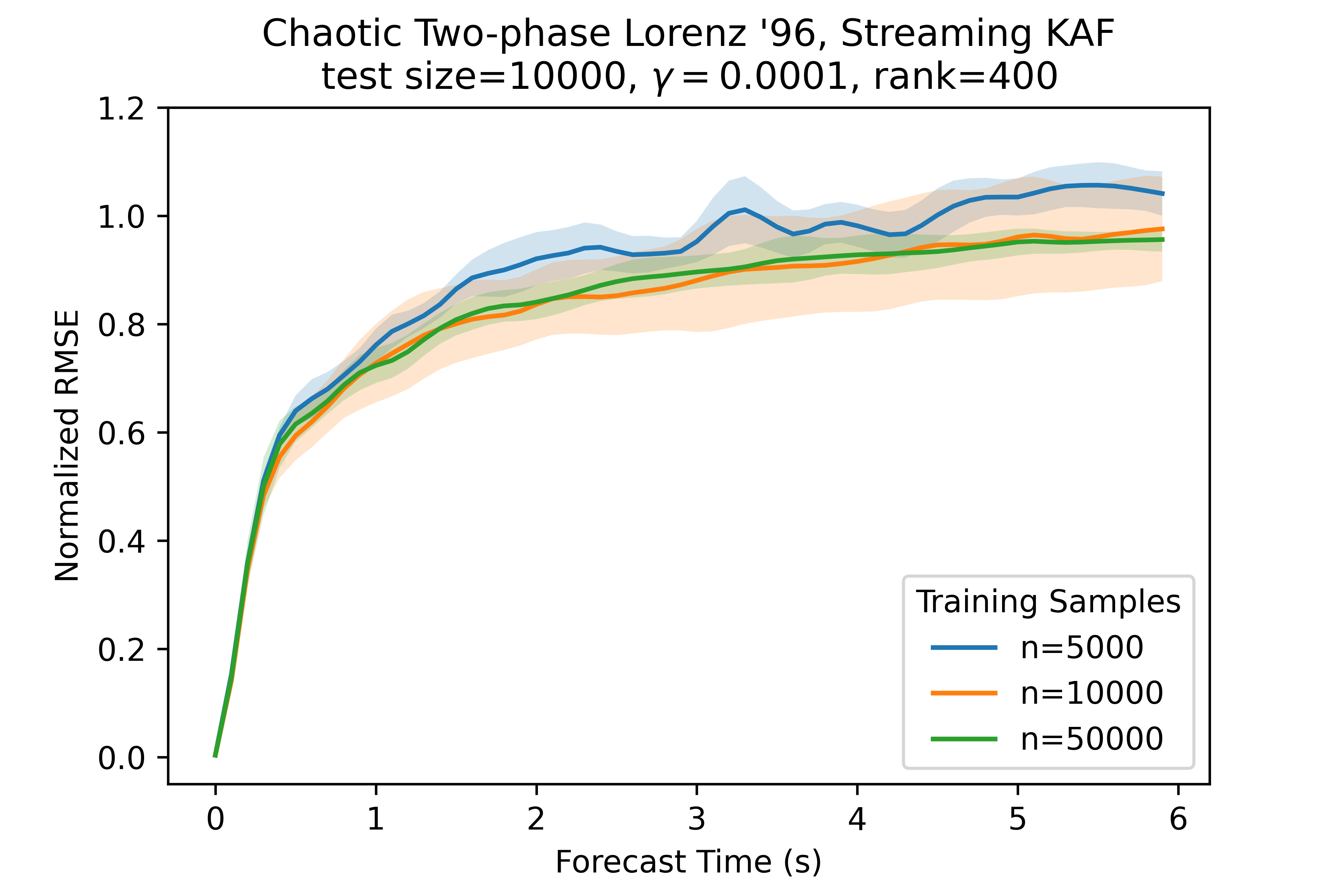}
    \caption{\textbf{Lorenz '96: Forecasting error versus amount of training data.}  Via streaming KAF, the average normalized RMSE for forecasting the first slow variable of periodic L96 [top], quasi-periodic L96 [bottom left], and chaotic L96 [bottom right] as a function of the number $n$ of training points.
    The regression model has dimension $\ell = 400$, and the number of features $s = \sqrt{n} \log(n)$.  The kernel inverse bandwidth $\gamma = 0.0001$ in the periodic and chaotic cases, while $\gamma = 0.01$ in the quasi-periodic case.  See~\cref{sec:l96}.}
    \label{fig:l96-rmse-vs-n}
\end{figure}

For the periodic regime ($F = 5$), forecasting is quite easy.  \Cref{fig:l96-rmse-vs-n} illustrates
the performance of streaming KAF as a function of the number $n$ of training samples.
The success of the method hardly varies as we increase $n$ from $5,000$ to $20,000$,
and the RMSE remains quite small over long time scales.

\begin{figure}[t]
    \centering
  \includegraphics[width=.45\textwidth]{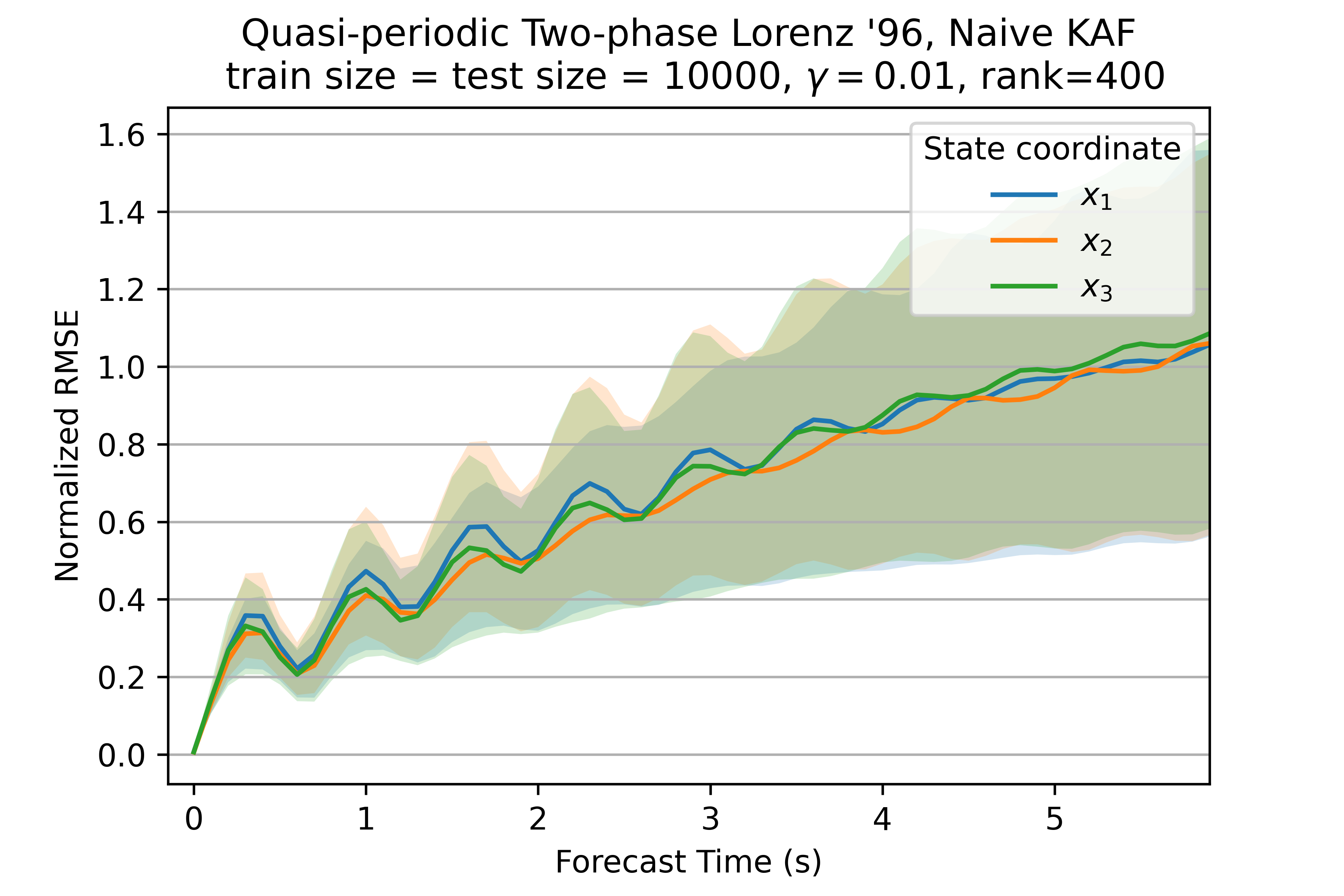} \hspace{1pc}
    \includegraphics[width=.45\textwidth]{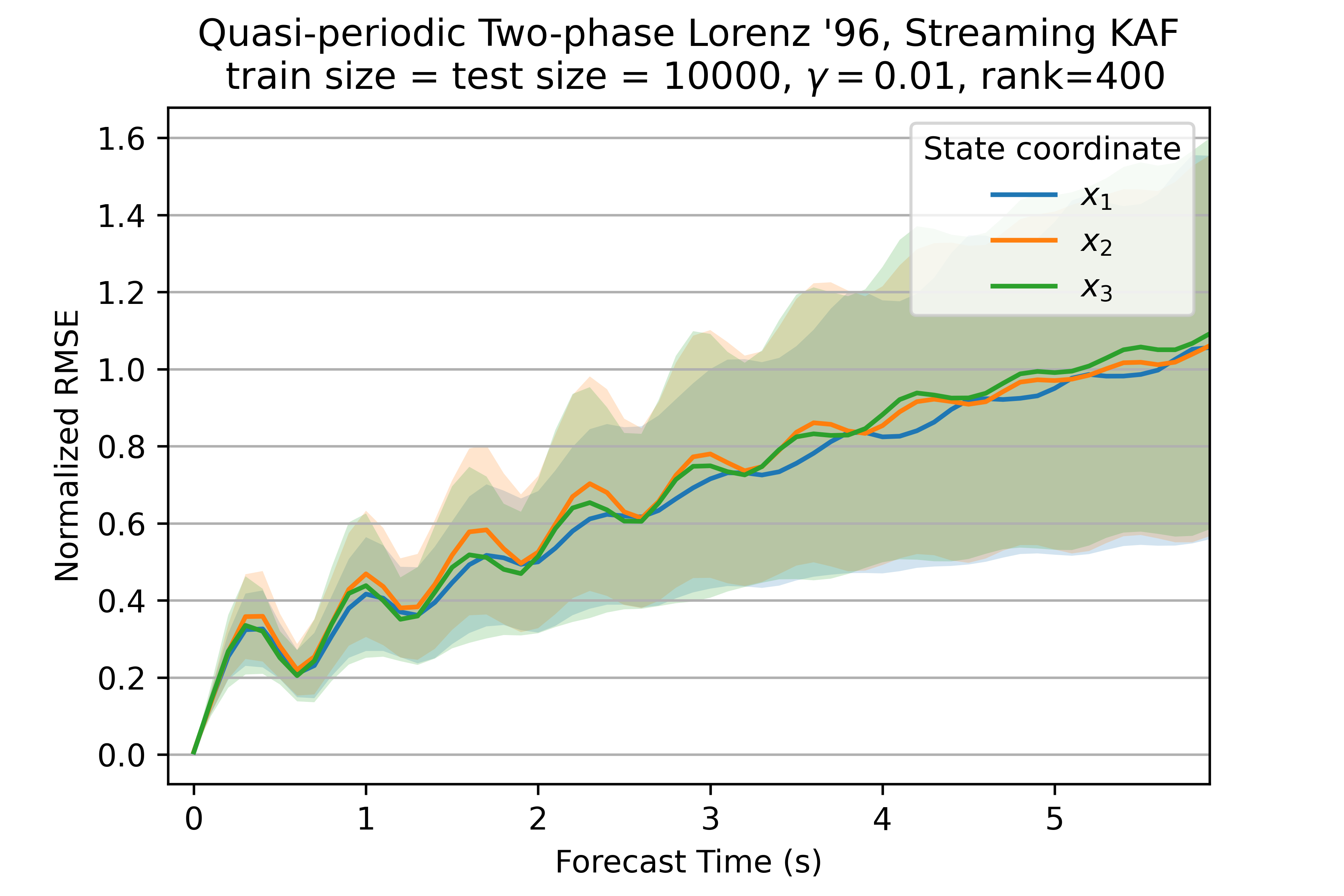}
    \caption{\textbf{Quasi-periodic Lorenz '96: Forecasting three slow variables.}  Average normalized RMSE for forecasting three slow coordinates $(x_1,x_2,x_3) = (\text{blue, orange, green})$ of quasi-periodic L96 via na{\"i}ve KAF [left] and streaming KAF [right] with $n = 10,000$ training points.
    The regression model has dimension $\ell = 400$, the kernel inverse bandwidth $\gamma = .0001$, and the number of features $s = \sqrt{n} \log(n)$.  See~\cref{sec:l96}.}
    \label{fig:l96-rmse-all-quasi}
\end{figure}

For the quasi-periodic regime ($F = 6.9$), the forecasting problem becomes more challenging.
For $n = 10,000$ training samples, 
\cref{fig:l96-rmse-all-quasi} shows that the na{\"i}ve and streaming methods have
similar forecasting performance for the first three slow variables.  
As we anticipate, the RMSE increases gradually with time.
\Cref{fig:l96-rmse-vs-n} displays the performance of streaming KAF as a function of the
number $n$ of training samples.  In this case, an increase in the number of samples
from $n = 10,000$ to $n = 50,000$ improves the performance moderately.
Note that the na{\"i}ve approach cannot benefit from the larger training set
because it does not scale to input of this size.

\begin{figure}[t]
    \centering
  \includegraphics[width=.45\textwidth]{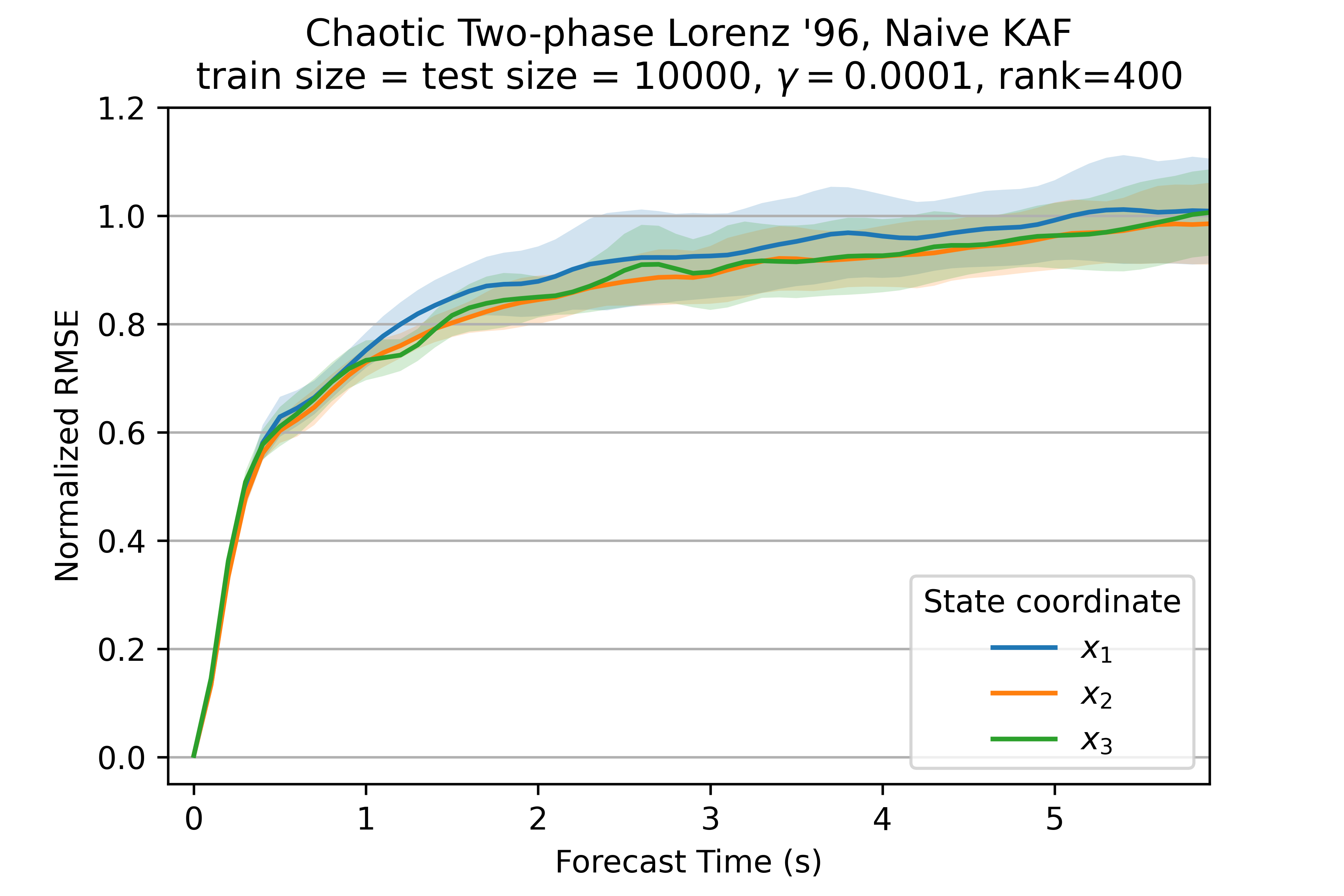} \hspace{1pc}
    \includegraphics[width=.45\textwidth]{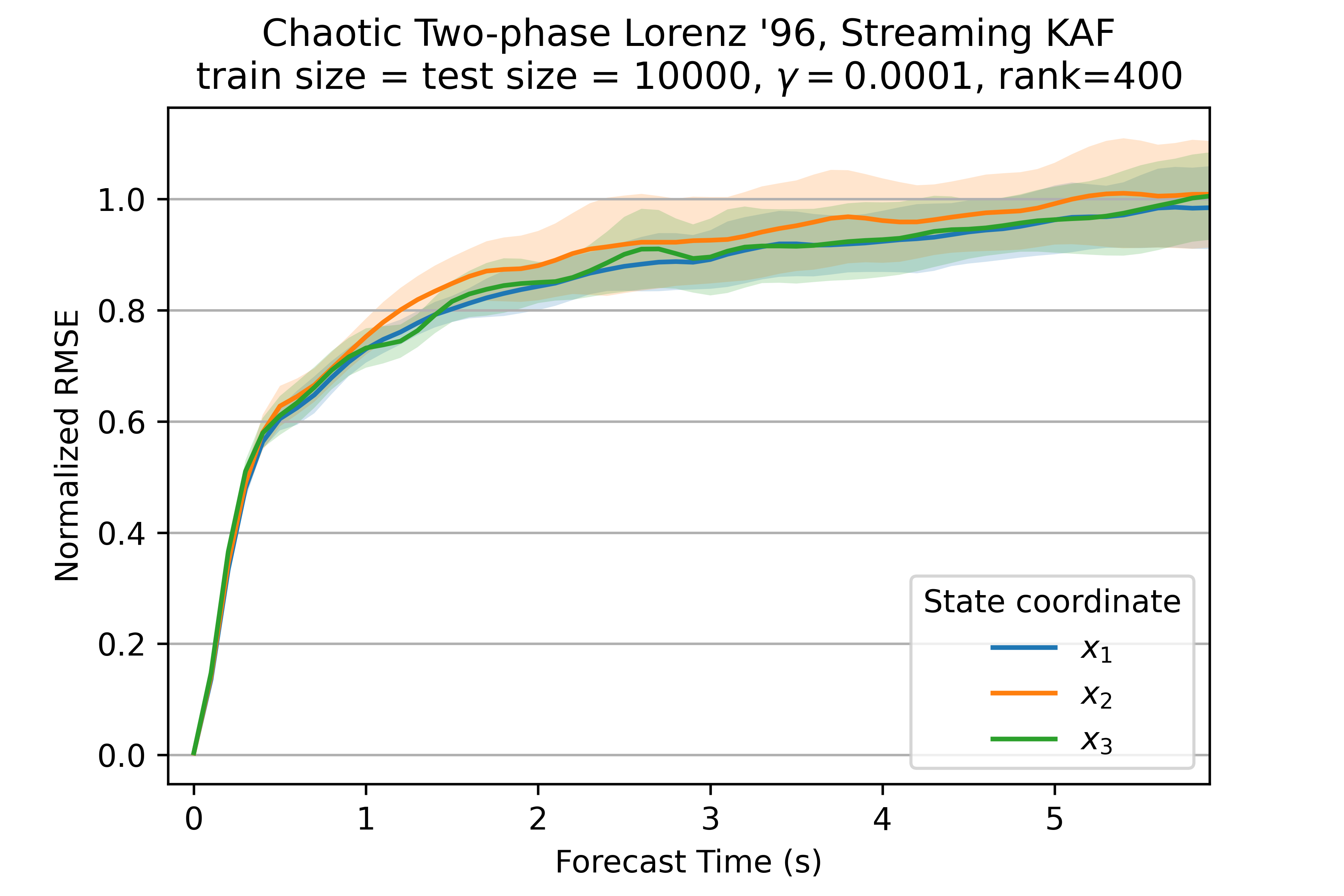}
        \caption{\textbf{Chaotic Lorenz '96: Forecasting three slow variables.}  Average normalized RMSE for forecasting three slow variables $(x_1,x_2,x_3) = (\text{blue, orange, green})$ of chaotic L96 via na{\"i}ve KAF [left] and streaming KAF [right] with $n = 10,000$ training points.
    The regression model has dimension $\ell = 400$, the kernel inverse bandwidth $\gamma = 0.0001$, and the number of features $s = \sqrt{n} \log(n)$.  See~\cref{sec:l96}.}
    \label{fig:l96-rmse-all-chaos}
\end{figure}

Last, we consider the chaotic regime ($F = 10$), where the forecasting problem is hard.
\Cref{fig:l96-rmse-all-chaos} indicates the na{\"i}ve and streaming methods
produce comparable forecasting results.  In both cases, the RMSE increases quite
quickly.  \Cref{fig:l96-rmse-vs-n} shows that streaming KAF can build models
from an increasing number $n$ of training samples, and it can attain an advantage
from the larger training set.

We conclude that streaming KAF and na{\"i}ve KAF have similar forecasting skill
in all three regimes, even though the streaming method makes several approximations.
At the same time, streaming KAF is far more economical, so it can exploit larger
sets of training data and thereby construct more accurate models.

\subsection{Hyperparameter specifications and sensitivity}

The streaming KAF method involves several hyperparameters: the kernel inverse bandwidth $\gamma$,
the dimension $\ell$ of the regression model, and the number $s$ of random features.
We performed a collection of experiments with the L63 and L96 data to gauge how much the hyperparameters affect
the quality of forecasts.

\subsubsection{Kernel bandwidth}
\label{sec:bandwidth-param}

The inverse bandwidth parameter $\gamma$ of the Gaussian RBF kernel is a notorious hyperparameter
that can have a significant impact on the performance of kernel methods.  One basic methodology
for selecting the bandwidth is the \textit{median rule}~\cite{2018Garreau}, which sets $\gamma^{-1/2}$ to be
the median pairwise distance among elements of a subsample from the dataset.
Other quantiles of the pairwise distance, such as the 0.1 and 0.9 quantiles, are sometimes employed. A different approach for bandwidth tuning~\cite{CoifmanEtAl08} leverages scaling relationships between the element sum of the $n \times n $ kernel matrix $\mtx K_{x,x}$ and $\gamma$.    

\begin{figure}[t]
    \centering
  \includegraphics[width=.4\textwidth]{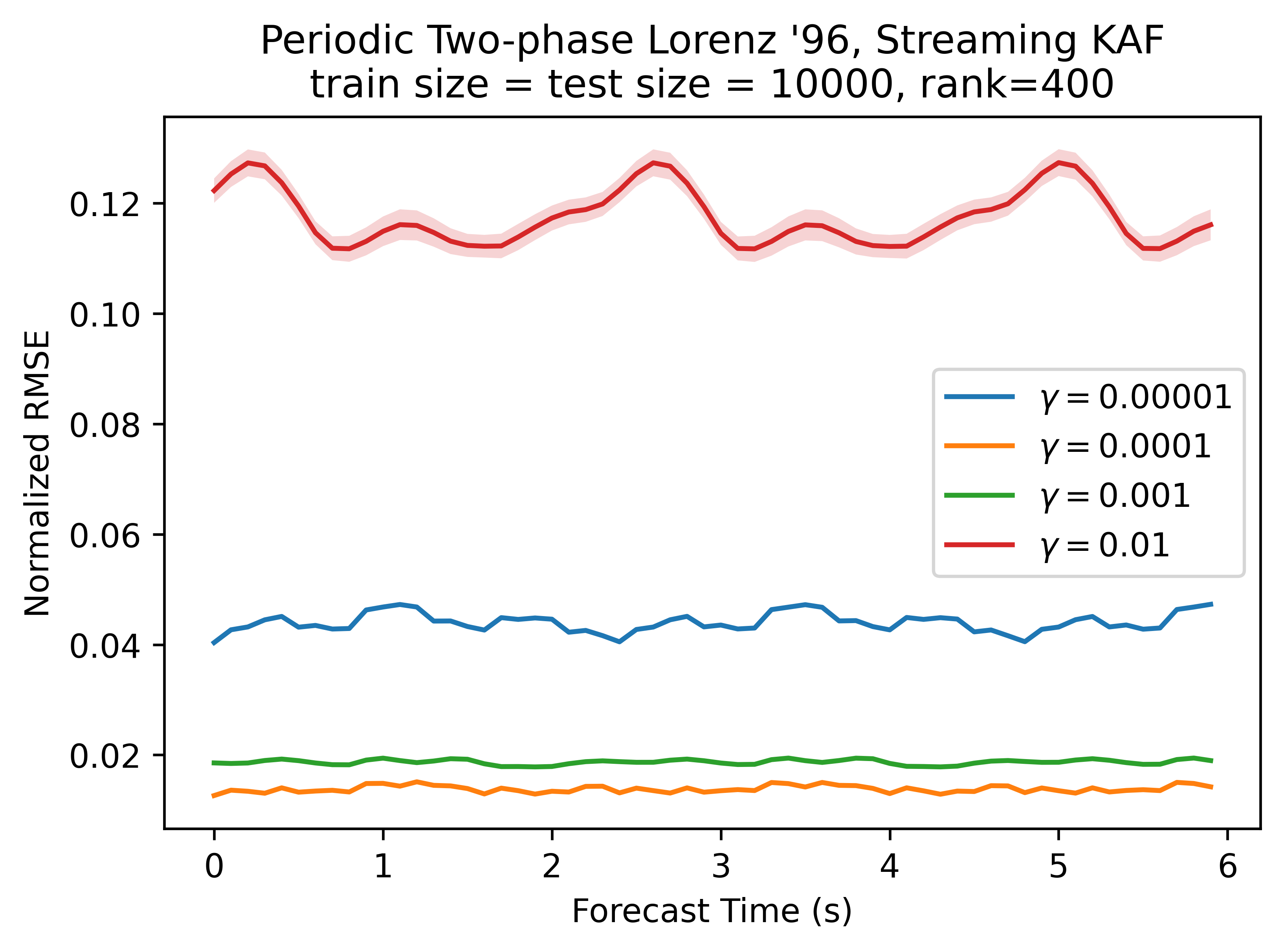} \hspace{1pc}
    \includegraphics[width=.45\textwidth]{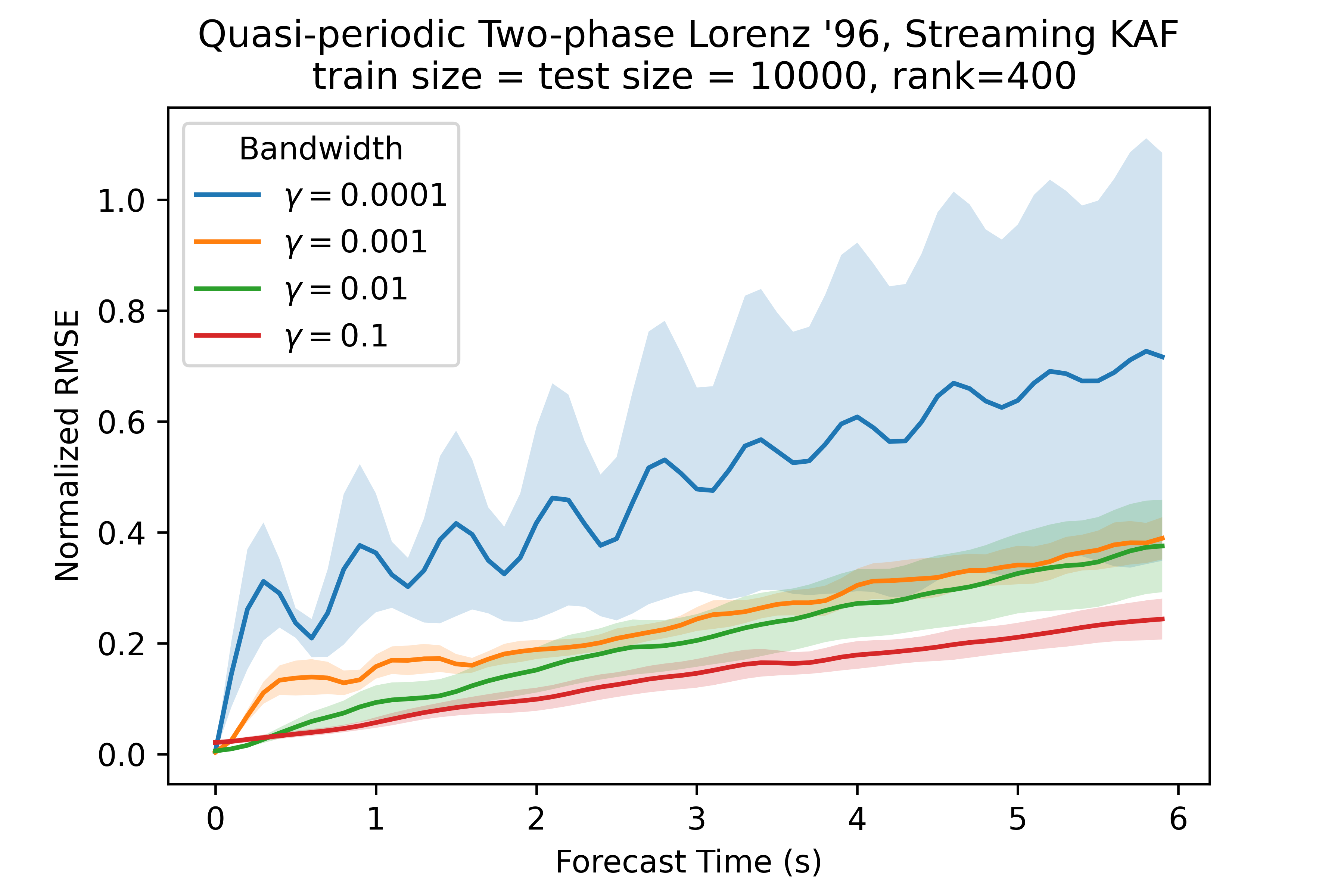}
    \caption{\textbf{Robustness to inverse bandwidth parameter:} For the L96 system in the periodic regime [left]
    and the quasi-periodic regime [right], the quality of the forecast is robust to the choice of
    $\gamma$.  In each case, we use $n = 10,000$ training samples, and the regression
    model has dimension $\ell = 400$.  See~\cref{sec:bandwidth-param}.}
    \label{fig:l96-gamma}
\end{figure}

In our experience, the KAF methodology is robust to the choice of inverse bandwidth parameter
in all problem regimes.  Indeed, the forecasting performance is similar over several orders of magnitude,
but tuning can have a modest effect.  See~\cref{fig:l96-gamma} for an illustration.  To obtain
better models from large training data, we invoke scaling laws for the bandwidth.
  
\subsubsection{Dimension of regression model}
\label{sec:dim-param}

To implement streaming KAF, we must choose the dimension, or rank, $\ell$ of the regression model.
When $\ell$ is too small, the model does not capture all of the dynamics.
Meanwhile, when $\ell$ is too large, we can introduce noise dimensions or
encounter numerical problems.
In this section, we outline some strategies for this task, and we will show
that the forecasting methodology is robust to the choice of this parameter.

One principled approach is to form the full covariance matrix $\mtx{C}_{xx} \in \RR^{s \times s}$
or $\mtx{C}_{uu} \in \RR^{s \times s}$ of the covariate data.
In this case, we can explicitly compute the eigenvalues
$(\lambda_1, \lambda_2, \dots, \lambda_s)$ of the matrix.
Then, we choose the truncation level $\ell$ so that we
capture, say, $99.9\%$ of the spectral content:
\begin{equation} \label{eqn:ell-rule}
\ell = \min \left\{ k \in \mathbb{N} : \sum\nolimits_{i=1}^k \lambda_i \geq 0.999 \cdot \sum\nolimits_{i=1}^s \lambda_i \right\}.
\end{equation}
This method is effective for a range of problems.  At the same time,
it imposes additional computational costs,
and it is not compatible with the streaming algorithm.

Instead, we typically prescribe the dimension $\ell$ of the regression model
in advance using prior knowledge about the problem or to work within our computational budget.
For example, in our medium-scale experiments, we make the choice $\ell = 400$,
which captures over $99.9\%$ of the spectral content of the computed covariance matrices.
Since we have included the ridge regularization $\mu \Id$ in the forecasting function,
we can insulate the algorithm from the negative impact of outsize $\ell$.

Given a conservative (i.e., large) initial value of $\ell$, randomized Nystr{\"o}m
produces an estimate for the first $\ell$ eigenvalues of the covariance matrix.
Using this estimate, we can apply the rule~\cref{eqn:ell-rule} a posteriori
to further reduce the dimension of the regression model.  This is often a good
compromise, but further research on principled methods would be valuable.

\begin{figure}[t]
    \centering
    \includegraphics[width=.45\textwidth]{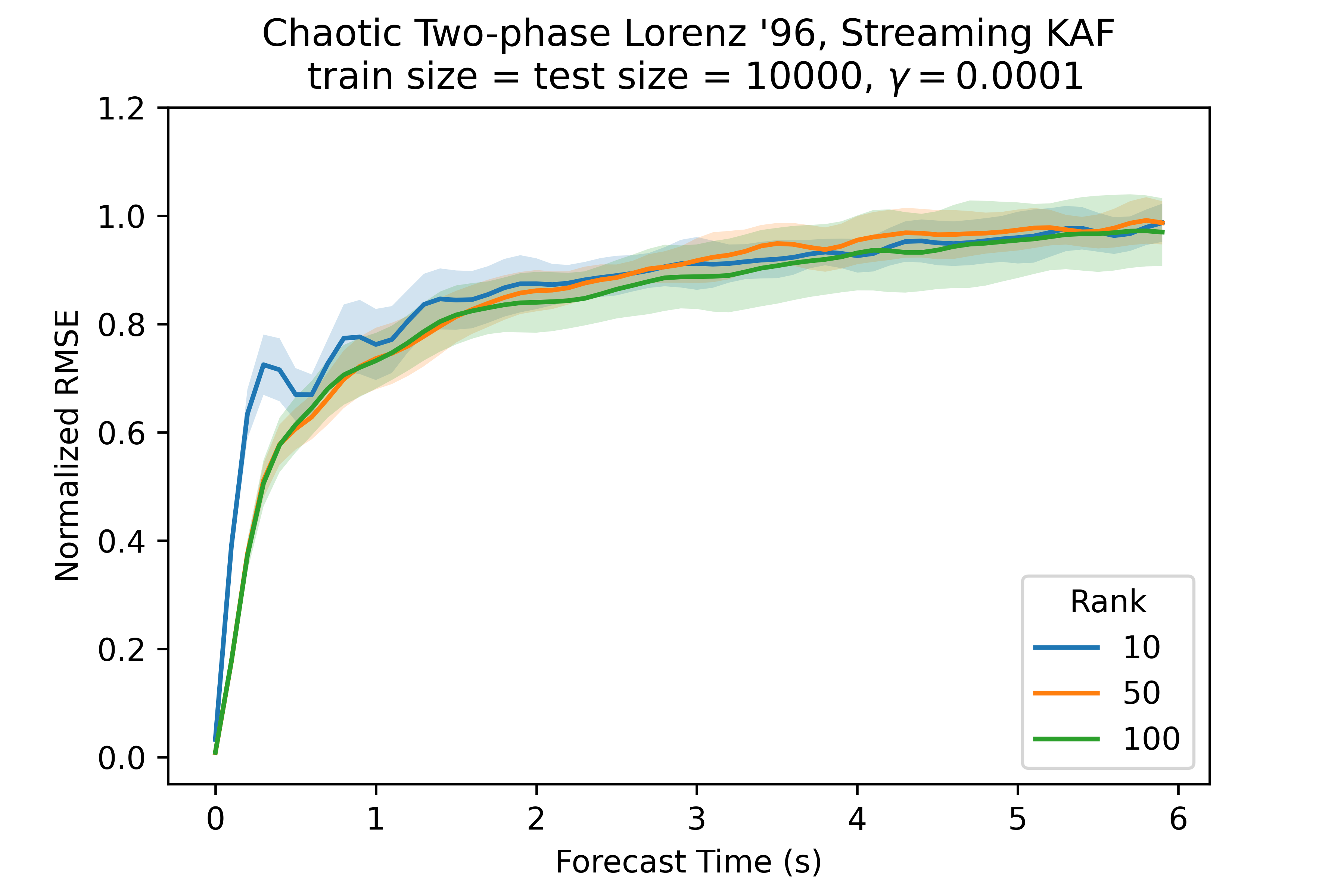} \hspace{1pc}
    \includegraphics[width=.45\textwidth]{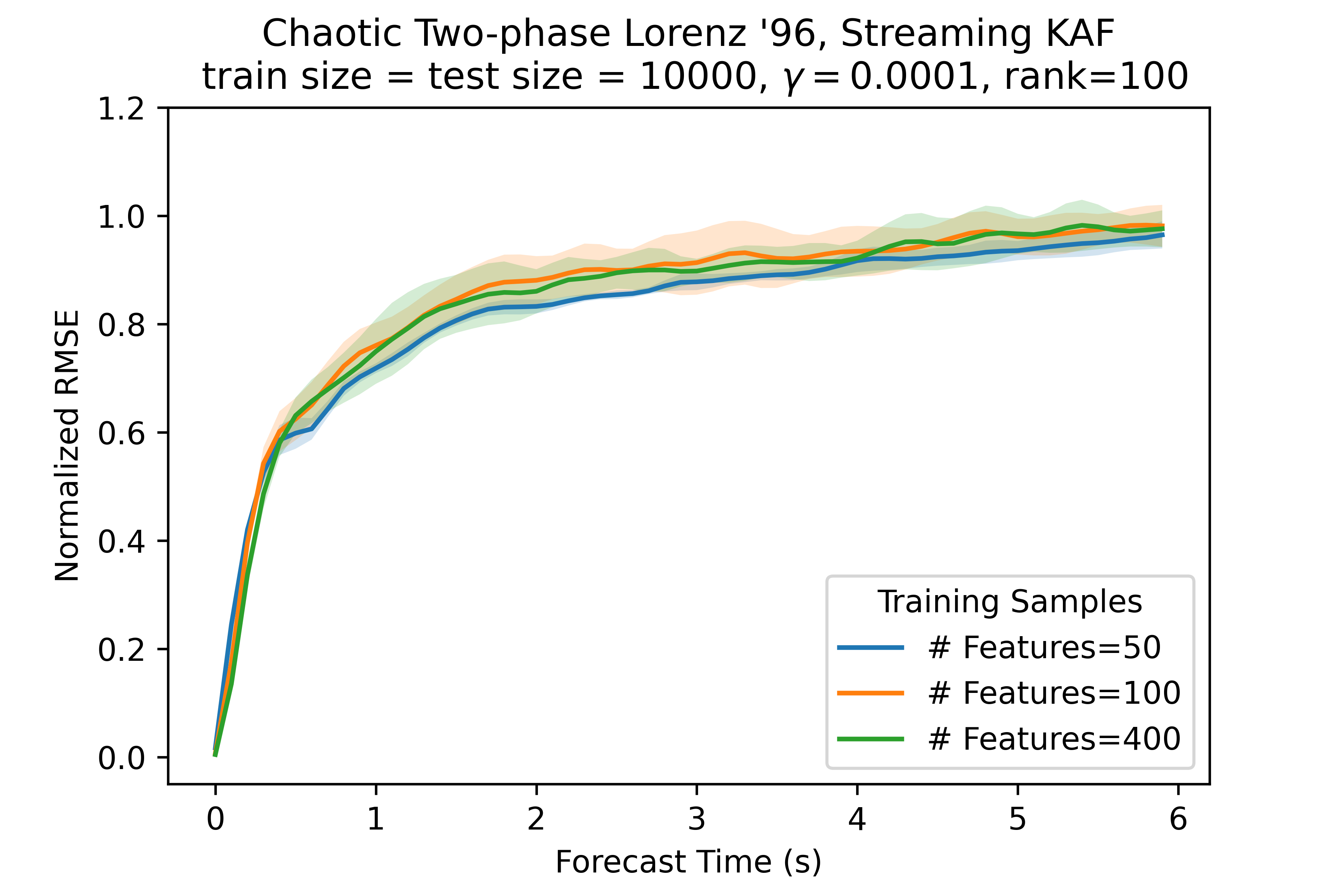}
    \caption{\textbf{Robustness to dimension of regression model and number of random features.} For the L96 system in the
    chaotic regime, the quality of the forecast is robust to the dimension $\ell$ of the regression model [left]
    and to the number $s$ of random features [right].    In the left panel, $s = 100$.  In the right panel, $\ell = 100$.
    In each case, we use $n = 10,000$ training samples, and the kernel inverse bandwidth $\gamma = 0.0001$.  See~\cref{sec:dim-param,sec:RFF}.}
    \label{fig:regression-features}
\end{figure}

Regardless, our numerical experiments indicate that streaming KAF forecast
is somewhat insensitive to the dimension $\ell$ of the regression model.
See~\cref{fig:regression-features} for some evidence.  For large training
data sets, we scale up the dimension $\ell$ to obtain more accurate forecasts.

\subsubsection{Number of random features}\label{sec:RFF}

The last parameter in the streaming KAF algorithm is the number $s$ of random features
that we use to approximate the kernel function.
As discussed in \cref{sec:rff}, the choice $s = \sqrt{n} \log(n)$ is theoretically
justified for kernel regression in a statistical setting.  In the majority of our experiments,
we adopt the value $s = \sqrt{n} \log(n)$, and we have found that the streaming KAF method
always performs well.  Furthermore, taking a larger number of random features does not seem to offer
any further benefit, and taking fewer random features is not detrimental.
See \cref{fig:regression-features} for evidence.

In the streaming setting, we may not know the number $n$ of training points in advance
and we do not want the model size to depend on the amount of input data,
so the prescription $s = \sqrt{n} \log(n)$ might be unappealing.  Our computational
work supports the recommendation that \textit{the number $s$ of random features may
be a small integer multiple of the dimension $\ell$ of the regression model}.
It would be interesting to understand this phenomenon better from both an
empirical and a theoretical point of view.

\begin{table}[t]
\footnotesize
    \centering
    \caption{\textbf{L63: Timing costs and error in forecasting.}  This table reports the time cost (in seconds) required to construct and evaluate a forecasting model using {\bf Streaming and Na{\"i}ve KAF algorithms}, along with the average normalized RMSE of the resulting models.  The covariate is the $3$-dimensional state of the L63 system, and the response is the first state variable after $0.5$ time units.  The number $n$ of training samples varies, and the number of random features $s = \sqrt{n}\log(n)$.  Reported test time is for {\bf making all $m = 10,000$ forecasts.}  See~\cref{sec:timings} for more details. }
    \label{tab:L63_train_time}
    \setlength\tabcolsep{1pt}
    \begin{tabular}{| c | l |r|r|r|r|r|r|}
        \hline
        & \diagbox[width=7em]{\textbf{Method} }{\scriptsize{${\bf (n,\ell, \gamma) = }$}} & \scriptsize{$(1e4, 4e2, .09)$}  & \scriptsize{$(5e4, 8e2, .18)$} & \scriptsize{$(1e5, 12e2,.27)$}  & \scriptsize{$(5e5,16e2,.36)$} & \scriptsize{(1e6, 24e2, .54)}  & \scriptsize{(5e6, 32e2, .72)} \\ \hline\hline
      \multirow{2}{*}{\textbf{Train}}
      	& Streaming  & .565 &  6.739  &  29.025 &   387.265  &   1588.570  & 26390.902 \\
      	& Na{\"i}ve  &  58.684 &   --- & --- & --- & --- & --- \\ \hline
      \multirow{2}{*}{\textbf{Test}}
        & Streaming   & .047 & .163 &    .267 &    .723  &   .936 & 2.703  \\
        & Na{\"i}ve  & 15.311  & --- & --- & --- & --- & --- \\ \hline\hline
      \multirow{2}{*}{\textbf{ RMSE }}
        & Streaming  & .262 &.177 &    .170 &    .107  & .065 & .047  \\
        & Na{\"i}ve   &  .228  &  --- & --- & --- & ---  & --- \\ \hline
    \end{tabular}
\end{table}

\begin{table}[t]
\footnotesize
    \centering
    \caption{\textbf{L63: Timing costs and error in forecasting using fewer random features.}  This table reports the time cost (in seconds) required to construct and evaluate a forecasting model {\bf using Streaming KAF}.  The setup is the same as in \cref{tab:L63_train_time}, but with model parameters $(s, \ell, \gamma) = (3200, 3200, .72)$ fixed for all experiments.   See~\cref{sec:timings} for more details. }
    \label{tab:small_s}
    \begin{tabular}{|c|r|r|r|r|r|r|r|r|r|}
        \hline
    & $\bm{n = 1e4}$ &  \textbf{5e4}   &   \textbf{1e5}  &  \textbf{5e5} & \textbf{1e6} & \textbf{5e6} \\ \hline\hline
      \multirow{1}{*}{\textbf{Train}}
 &  35.101   & 43.051	 &  55.079 &  138.083 &  253.435 & 1061.187 \\ \hline
      \multirow{1}{*}{\textbf{Test}}
  &  .230  &.220  &    .245 &  .229 &    .221& .219  \\ \hline\hline
      \multirow{1}{*}{\textbf{ RMSE }}
  & .408  &   .125 &      .119  &   .086  &  .075  & .089  \\ \hline
    \end{tabular}
\end{table}

\subsection{Timing comparisons} \label{sec:timings}

We have demonstrated that streaming KAF constructs accurate forecasting models in a range of scenarios.
Therefore, we may turn our attention to the computational costs of training and forecasting.
\Cref{tab:L63_train_time} compares the runtimes of na{\"i}ve  KAF and streaming KAF,
and it charts the average normalized RMSE of the resulting models.

These experiments are based on the L63 data.  We forecast the first state variable
from the full set of three state variables.  The forecast horizon is fixed at $q = 0.5$ time units.
The number $n$ of training samples varies, while the number of test samples remains fixed at $m = 10,000$.  
The kernel inverse bandwidth $\gamma = 0.09$, and the dimension of the regression model
grows from $\ell = 400$ to $\ell = 3200$ in rough proportion to $\log n$.
For the streaming method, the number of random features $s = \sqrt{n} \log(n)$
also increases with the size of the training data.
We report the average RMSE over five test runs.

To be clear, the training time includes the full cost of computing the
weight matrix $\check{\mtx{W}}_{q, \ell}$ for a single real-valued response at a
single forecast horizon $q$.  This cost includes the evaluation of random features,
formation of the covariance matrices, the streaming PCA computation, and the matrix product.
For making a forecast, the timing reflects the full cost of computing
$m = 10,000$ real-valued responses for the fixed time horizon $q$,
including the evaluation of random features and the matrix product.

For small problems, we see that the training time for streaming KAF is 100--200$\times$
faster than na{\"i}ve KAF.  The test time for streaming KAF is 300--400$\times$ faster,
and the models achieve similar RMSE.  For large problems, na{\"i}ve KAF is unable
to produce a forecasting model.  Meanwhile, streaming KAF can build a forecasting model
from $n = 5 \cdot 10^6$ training samples in less than two hours on a laptop,
and this model can produce a single real-valued forecast in about $0.0003\mathrm{s}$.
As the amount of training data increases, the RMSE of the forecasting models
continues to improve, which underscores how important it is to develop a
scalable algorithm.

Out of a sense of fair play, we used the theoretically supported number
$s = \sqrt{n} \log(n)$ of random features.
If we adopt our empirical recommendation $s = \mathrm{Const} \cdot \ell$,
the timings improve markedly without sacrificing much accuracy.
\Cref{tab:small_s} displays the runtimes and average normalized RMSE for streaming KAF 
under the same experimental set-up as in \Cref{tab:L63_train_time}, but with the regression
dimension and number of random features fixed at $(\ell, s) = (3200, 3200)$.  The user
may judge whether the speedup warrants the modest loss in RMSE.


\section{Comparison with related work}
\label{sec:related}

Several other techniques for data-driven {prediction} have been proposed and studied recently. Here, we comment on the mathematical and computational characteristics of these approaches in relation to streaming KAF, focusing on methods that employ aspects of linear operator theory {or randomized linear algebra}. Within this context, forecasting techniques can be broadly classified as reduced modeling approaches (i.e., methods that construct a surrogate dynamical system from observed data) and regression approaches (i.e., supervised learning techniques for estimating covariate--response relationships).

\subsection{{Forecasting methodologies}}

Examples of reduced modeling techniques are linear inverse models \cite{Penland89}, (extended) DMD \cite{RowleyEtAl09,Schmid10,WilliamsEtAl15}, and methods for approximating the Koopman generator \cite{Giannakis19,DasEtAl21,KlusEtAl20b}. These methods formally assume that the training data have a (deterministic) Markovian evolution. That assumption is clearly satisfied under the autonomous dynamics in~\eqref{eqn:state-trajectory} if the training data are snapshots $\mtx x_0, \mtx x_1, \ldots$ of the full system state in $\mathbb R^d$. On the other hand,  if we have access to samples $u(\mtx x_0), u(\mtx x_1), \ldots$ of a covariate $ u : \mathbb R^d \to \mathbb R^{d'}$ with $d' < d$, then training data are generally non-Markovian (unless $u$ happens to lie in a Koopman-invariant subspace).  Approaches for overcoming non-Markovianity include dimension augmentation through delay-coordinate maps \cite{Sauer93,BruntonEtAl17,GilaniEtAl21} and incorporation of memory terms using the Mori-Zwanzig formalism \cite{GouasmiEtal17,GutierrezEtAl21}.

Other approaches model the observed data as realizations of a stochastic process. For example, techniques based on Ulam's method~\cite{DellnitzJunge99,JungeKoltai09} estimate the transfer operator of a dynamical system (which is a dual operator to the Koopman operator, acting on probability measures) in a basis of indicator functions associated with a partition of state space. The diffusion forecasting technique \cite{BerryEtAl15} estimates the evolution semigroup associated with a stochastic differential equation (SDE) on a manifold in a smooth data-driven basis of kernel eigenfunctions learned through the diffusion maps algorithm \cite{CoifmanLafon06}. Extensions of DMD to random dynamical systems \cite{CrnjaricEtAl19} and SDEs \cite{ArbabiSapsis21} have also been proposed recently.

A common aspect of reduced modeling techniques is that they learn a surrogate model of the dynamics from time series data. Often, in order to make a forecast to a horizon of $q$ time units, these models are trained on a shorter timestep $q'< q $ and iteratively applied $q/q'$ times to reach the desired horizon. This approach is attractive because it allows simulation of the long-term statistical behavior of the system (assuming that the training phase was successful). 

In contrast, regression-based methodologies usually operate by constructing a forecast function at a \emph{fixed} lead time (or a family of independent forecast functions up to a desired lead time), and they evaluate the forecast once on the initial data to yield a prediction. This approach offers greater generality than reduced modeling approaches, since Markovianity of the covariate--response observables is not required, nor is it required that the covariate and response lie in a Koopman-invariant subspace \cite{GilaniEtAl21}. Indeed, as discussed in \cref{sec:koop}, KAF yields asymptotically optimal predictions (in the $L_2$ or RMSE sense) in the large-data limit in the form of the conditional expectation of the Koopman-evolved response conditioned on the covariate. Yet, at the same time, the conditional expectation may not be a good approximation for actual dynamical trajectories, which makes direct regression approaches unsuitable for simulating the statistical behavior of the system (despite yielding RMSE-optimal forecasts). For further details, see the paper \cite{BurovEtAl21}, which studies applications of KAF to multiscale systems with averaging and homogenization limits.  A recent paper \cite{HamziOwhadi21} has explored applications of kernel learning \cite{OwhadiYoo19} to forecasting with kernel regression. 

All of the above approaches are purely data-driven, in the sense that they only use time-ordered data snapshots as inputs, without requiring knowledge of the equations of motion.  Yet, in many applications, full or partial knowledge of the equations of motion \emph{is} available, and it is natural to design methods that take advantage of that knowledge \cite{HamiltonEtAl15}. An example is the ``lift and learn'' framework~\cite{qian2020lift} which employs a mapping to transport the data to a higher-dimensional space where the system is quadratic. Unlike the Koopman operator, the existence of a finite-dimensional quadratic representation of the system dynamics is not universally guaranteed, but can be constructed for many systems encountered in physical and engineering applications \cite{gu2011qlmor} if the equations of motion are known. The approach of \cite{qian2020lift} leverages the quadratic structure of the system in the lifted space by employing a projection that is compatible with quadratic nonlinearities (see also \cite{peherstorfer2016data}). In this manner, the reduced model is compatible with the ``physics'' of the lifted model. In \cite{qian2020lift}, the projection is obtained from the \emph{proper orthogonal decomposition} (POD) \cite{HolmesEtAl96}, which computes a low-rank approximation to the \emph{autocorrelation} matrix $\mtx X \mtx X^{\tsp} \in \mathbb R^{d\times d}$ rather than the covariance matrix $\mtx X^{\tsp} \mtx X \in \mathbb R^{n\times n} $.  The randomized singular value decomposition is  also used within this forecasting framework to build a scalable implementation \cite{mcquarrie2021data}.

{Note that the eigenvectors of $\mtx X \mtx X^{\tsp}$ are spatial vectors in $\mathbb R^d$. In DMD, the analogous objects are the eigenvectors of the matrix $\mtx A$ in~\eqref{eqDMDA}, called Koopman modes \cite{RowleyEtAl09}, which can also be employed for model reduction. The KAF approach can be thought of as being ``dual'' to these methods in that it employs $n\times n$ kernel matrices which are discretizations of operators acting on spaces of observables of the system (rather than spatial patterns in $\mathbb R^d$).}

\subsection{Streaming algorithms for kernel computation}

The machine learning literature contains a substantial body of work on kernel methods, techniques for combining kernels with random features, and methods for implementing these algorithms in a streaming setting.  This space is not adequate for a comprehensive summary of this vast field.  We recommend the book~\cite{scholkopf2002learning} as a foundational reference on kernel methods in machine learning.

The RFF technique~\cite{2008Rahimi} was developed to accelerate kernel computations. 
There are a substantial number of papers that use RFF for KRR, such as~\cite{AKM+17:Random-Fourier, 2017Rudi}, but we are not aware of a paper that uses random features for \textit{streaming} kernel regression.

There are also several papers that combine RFF with streaming PCA algorithms to obtain streaming KPCA algorithms.  In particular, Ghashami et al.~\cite{2016Ghashami} apply the frequent directions method~\cite{GLPW16:Frequent-Directions}, while Ullah et al.~\cite{2018Ullah} use Oja's algorithm~\cite{Oja82:Simplified-Neuron}.  Henriksen \& Ward~\cite{2019Henriksen} have developed an adaptive extension of Oja's algorithm that is significantly more robust.  Tropp and coauthors have proposed to use the randomized Nystr{\"o}m method for streaming PCA~\cite{2017Tropp}, perhaps in combination with random features~\cite[Sec.~19.3.5]{martinsson2020randomized}.  Our numerical work suggests that the Nystr{\"o}m method is more accurate and more reliable than the alternatives in the context of streaming KAF.

We have also investigated the performance of streaming KAF using AdaOja~\cite{2019Henriksen} for the streaming PCA computation.  In our experience, this approach can be competitive, especially in cases where the spectrum of the covariance matrix decays slowly.  See~\cite{2021HenriksenThesis} for a detailed report.  

\section{Conclusions} 

Kernel analog forecasting is a regression-based approach to forecasting dynamical systems that offers a theoretical guarantee of asymptotically optimal predictions (in the $L_2$ or RMSE sense) in the large-data limit.  By incorporating two randomized approximation techniques from numerical linear algebra---random Fourier features and the randomized Nystr{\"o}m method---we developed a streaming implementation of kernel analog forecasting.  This approach makes it possible to build forecasting models from large data sets where the KAF methodology is theoretically justified.  Our experiments indicate that streaming KAF has the potential to unlock the promise of KAF
as a general data-driven, non-parametric tool for making predictions of dynamical systems.

\section*{Acknowledgments}
We are thankful for helpful feedback from Eliza O'Reilly and Ethan Epperly. DG thanks the department of Computing and Mathematical Sciences at the California Institute of Technology for hospitality during a sabbatical in 2017/18 where part of this work was initiated.    

\bibliographystyle{siamplain}
\bibliography{references,bibliography}

\end{document}

%% file: ex_shared.tex
\usepackage{lipsum}
\usepackage{amsfonts}
\usepackage{graphicx}
\usepackage{epstopdf}
\usepackage{comment}

\usepackage[noend]{algpseudocode}
\usepackage{algorithm,algorithmicx}

\algrenewcommand\alglinenumber[1]{\sf\scriptsize\color{blue}{#1}}
\algrenewcommand\algorithmicrequire{\textbf{Input:}}
\algrenewcommand\algorithmicensure{\textbf{Output:}}

\algdef{SE}[SUBALG]{Indent}{EndIndent}{}{\algorithmicend\ }%
\algtext*{Indent}
\algtext*{EndIndent}

\ifpdf
  \DeclareGraphicsExtensions{.eps,.pdf,.png,.jpg}
\else
  \DeclareGraphicsExtensions{.eps}
\fi

\usepackage{enumitem}
\usepackage{fancyhdr}
\setlist[enumerate]{leftmargin=.5in}
\setlist[itemize]{leftmargin=.5in}

\newsiamremark{remark}{Remark}
\newsiamremark{hypothesis}{Hypothesis}
\crefname{hypothesis}{Hypothesis}{Hypotheses}
\newsiamthm{claim}{Claim}

\newsiamthm{prop}{Proposition}
\newsiamthm{defn}{Definition}

\headers{Streaming Kernel Analog Forecasting}{D.~Giannakis, A.~Henriksen, J.~A.~Tropp, R.~Ward}

\title{Learning to Forecast Dynamical Systems from Streaming Data\thanks{Submitted to the editors DATE.
\funding{DG acknowledges support from NSF DMS 1854383 and ONR MURI N00014-19-1-242. AH was funded by AFOSR MURI FA9550-19-1-0005, NSF DMS 1952735.  JAT was supported by ONR N00014-18-1-2363 and NSF DMS 1952777.  RW acknowledges support from AFOSR MURI FA9550-19-1-0005, NSF DMS 1952735, and NSF IFML 2019844.}}}

\author{Dimitrios Giannakis\thanks{Dartmouth College, Hanover, NH (\email{dimitrios.giannakis@dartmouth.edu}); Courant Institute of Mathematical Sciences, New York University, New York, NY.} \and
Amelia Henriksen\thanks{Oden Institute, University of Texas at Austin, Austin, TX (\email{amelia@oden.utexas.edu}).}
\and Joel A.~Tropp\thanks{Computing and Mathematical Sciences, California Institute of Technology, Pasadena, CA (\email{jtropp@cms.caltech.edu}).}
\and Rachel Ward\thanks{Department of Mathematics, University of Texas at Austin, Austin, TX  (\email{rward@math.utexas.edu}).}
}

\usepackage{amsopn}
\DeclareMathOperator{\diag}{diag}

\usepackage{mathtools}